\pdfoutput=1
\documentclass[final,onefignum,onetabnum]{siamart220329}

\usepackage{amssymb}
\usepackage{booktabs}
\usepackage{array}
\usepackage{algorithm}
\usepackage{algpseudocode}
\usepackage{placeins}
\usepackage{microtype}

\numberwithin{algorithm}{section}
\crefname{algorithm}{Algorithm}{Algorithms}
\Crefname{algorithm}{Algorithm}{Algorithms}

\newsiamthm{assumption}{Assumption}
\newsiamremark{remark}{Remark}

\makeatletter
\AtBeginDocument{%
  \@ifundefined{MakeLinkTarget}{}{%
    \def\cref@old@refstepcounter#1{%
      \stepcounter{#1}%
      \MakeLinkTarget{#1}%
      \protected@edef\@currentlabel{\csname p@#1\endcsname\csname the#1\endcsname}}%
    \def\refstepcounter@noarg#1{%
      \cref@old@refstepcounter{#1}%
      \cref@constructprefix{#1}{\cref@result}%
      \@ifundefined{cref@#1@alias}%
        {\def\@tempa{#1}}%
        {\def\@tempa{\csname cref@#1@alias\endcsname}}%
      \protected@edef\cref@currentlabel{%
        [\@tempa][\arabic{#1}][\cref@result]%
        \csname p@#1\endcsname\csname the#1\endcsname}}%
    \def\refstepcounter@optarg[#1]#2{%
      \cref@old@refstepcounter{#2}%
      \cref@constructprefix{#2}{\cref@result}%
      \@ifundefined{cref@#1@alias}%
        {\def\@tempa{#1}}%
        {\def\@tempa{\csname cref@#1@alias\endcsname}}%
      \protected@edef\cref@currentlabel{%
        [\@tempa][\arabic{#2}][\cref@result]%
        \csname p@#2\endcsname\csname the#2\endcsname}}%
  }%
}
\makeatother

\newcommand{\curl}{\nabla\times}
\newcommand{\divg}{\nabla\cdot}
\newcommand{\grad}{\nabla}
\newcommand{\Hone}{{H^1}}
\newcommand{\Hcurl}{H(\mathrm{curl})}
\newcommand{\Hdiv}{H(\mathrm{div})}
\newcommand{\RR}{\mathbb{R}}
\newcommand{\bfu}{\mathbf{u}}
\newcommand{\bfv}{\mathbf{v}}
\newcommand{\bfw}{\mathbf{w}}
\newcommand{\bff}{\mathbf{f}}
\newcommand{\bfn}{\mathbf{n}}
\newcommand{\bfx}{\mathbf{x}}
\newcommand{\calT}{\mathcal{T}}
\newcommand{\calE}{\mathcal{E}}
\newcommand{\ND}{\mathrm{ND}}
\newcommand{\RT}{\mathrm{RT}}
\newcommand{\true}{\mathrm{true}}
\newcommand{\loc}{\mathrm{loc}}
\newcommand{\ho}{\mathrm{ho}}
\newcommand{\lorr}{\mathrm{lor}}

\headers{MATRIX-FREE LOR FOR VARIABLE-ORDER DE RHAM COMPLEX}{C. QIN}

\title{Matrix-Free Low-Order-Refined Preconditioning for Variable-Order Nonconforming de~Rham Discretizations}

\author{Ce Qin\thanks{School of Physics and Electronic Information Engineering, Henan Polytechnic University, Jiaozuo 454000, China (\email{ce.qin@hpu.edu.cn}).}}

\begin{document}
\maketitle

\begin{abstract}
High-order $hp$-adaptive discretizations of the three-dimensional $\Hone$--$\Hcurl$--$\Hdiv$ de~Rham complex deliver high accuracy per unknown on locally refined hexahedral meshes, yet assembled matrices and high-order algebraic preconditioners incur $\mathcal{O}(p_K^6)$ element-local memory and setup costs. Building on variational-restriction low-order-refined (LOR) preconditioning, we develop a parallel matrix-free solver framework for variable-order nonconforming de~Rham discretizations built around three components: three-dimensional mixed-order hanging-interface constraint completion, explicit sparse auxiliary-space interpolants preserving signed constraint weights, and order-grouped partial assembly coupled to heterogeneous-degree LOR solvers. Iterative minimum-rule propagation and parallel internal-edge compatibility elimination build globally conforming discrete subcomplexes without over-constraining master faces. Grouping elements by polynomial degree enables sum-factorized operator evaluation with $\sum_{K\in\calT_h^\loc}\mathcal{O}(p_K^3)$ element-local storage and $\sum_{K\in\calT_h^\loc}\mathcal{O}(p_K^4)$ element-local work. On each element's Gauss--Lobatto sub-grid, element-local $\mathcal{O}(p_K^3)$ direct-sparse LOR assembly pairs topological differential matrices with explicit N\'ed\'elec and corner-quadrature Raviart--Thomas interpolants restricted to owner entities, avoiding the constant-field mismatch caused by coordinate-based reconstructions under negative high-order Lagrange weights. Under compatible geometric traces and a conforming subcomplex, we prove exact discrete gradient and curl commutativity, exact constant-flux integration on warped bilinear sub-faces, and uniform global spectral equivalence from element-wise bounds. Parallel experiments up to $p=6$, 1{,}024 cores, and 90 million unknowns demonstrate degree- and contrast-robust convergence, up to a $43\%$ reduction in $\Hdiv$ iterations over coordinate reconstruction on mixed-order meshes, and up to $13\times$ speedup over assembled solvers.
\end{abstract}



\section{Introduction}
\label{sec:intro}
Structure-preserving discretizations of the three-dimensional de~Rham sequence
\begin{equation}
\Hone(\Omega) \xrightarrow{\;\grad\;} \Hcurl(\Omega) \xrightarrow{\;\curl\;} \Hdiv(\Omega)
\label{eq:derham-cont}
\end{equation}
underpin finite element methods in heat transport, computational electromagnetics, and porous-media flow \cite{ArnoldFalkWinther2006,ArnoldFalkWinther2010,BochevHyman2006}. Matching continuous Lagrange, N\'ed\'elec \cite{Nedelec1980}, and Raviart--Thomas \cite{RaviartThomas1977} spaces preserves the discrete kernels of the curl and divergence, ensuring stability in mixed and vanishing-reaction ($\beta\to 0$) regimes on affine and curvilinear hexahedra \cite{ArnoldBoffiBonizzoni2015,BoffiBrezziFortin2013,Hiptmair2002,Monk2003}. We consider the symmetric positive definite (SPD) model problems
\begin{subequations}\label{eq:models}
\begin{align}
-\divg(\alpha \grad u) + \beta u &= f \quad \text{in } \Omega, & u &= u_D \text{ on } \partial\Omega, \label{eq:model-h1}\\
\curl(\alpha \curl \bfu) + \beta \bfu &= \bff \quad \text{in } \Omega, & \bfn \times \bfu &= \bfn \times \bfu_D \text{ on } \partial\Omega, \label{eq:model-hcurl}\\
-\grad(\alpha \divg \bfu) + \beta \bfu &= \bff \quad \text{in } \Omega, & \bfn \cdot \bfu &= \bfn \cdot \bfu_D \text{ on } \partial\Omega, \label{eq:model-hdiv}
\end{align}
\end{subequations}
with positive coefficients $\alpha,\beta$ that may jump by orders of magnitude across material interfaces.

When solutions exhibit corner singularities or localized layers, $hp$-adaptive discretizations combine nonconforming octree refinement with element-wise polynomial degrees $p_K$ to achieve rapid convergence \cite{BabuskaGuo1992,BabuskaSuri1994,Demkowicz2006,DemkowiczKurtzPardo2007,MelenkWohlmuth2001,Schwab1998}. Scaling such discretizations across \eqref{eq:derham-cont} faces three coupled bottlenecks. First, conformity across parallel $k$-irregular hanging interfaces requires algebraic prolongations that handle edge-orientation signs, distributed degree consensus, and coupled edge--face constraints \cite{BangerthBursteddeHeisterKronbichler2011,CervenyDobrevKolev2019,HoeflerSiebertLumsdaine2010NBX,KinnewigWickBeuchler2025}. Second, full matrix assembly incurs $\mathcal{O}(p_K^6)$ element-local storage per hexahedron. This motivates sum-factorized partial assembly \cite{AndersonEtAl2021MFEM,Brown2010,DevilleFischerMund2002,KarniadakisSherwin2005,KolevCEED2021,KronbichlerKormann2012,KronbichlerWall2018,Orszag1980Spectral}, where scalar and discontinuous Galerkin (DG) codes handle variable degrees by batching equal-degree elements \cite{CantwellSherwinKirbyKelly2011,MunchHeisterPrietoKronbichler2023}. Third, BoomerAMG \cite{HensonYang2002BoomerAMG,XuZikatanov2017} and Hiptmair--Xu auxiliary-space solvers (AMS and ADS in hypre \cite{ArnoldFalkWinther2000,FalgoutYang2002,Hiptmair1998,HiptmairXu2007,KolevVassilevski2009,KolevVassilevski2012ADS}) scale well at lowest order, whereas applying them directly to assembled high-order matrices causes severe fill-in and $p$-dependent degradation \cite{SundarStadlerBiros2015}; this has motivated $p$-multigrid, Schwarz, and patch smoothers \cite{BabuskaCraigMandelPitkaranta1991,BrubeckFarrell2022,Casarin1997,LottesFischer2005,Mandel1990,Pavarino1994,ToselliWidlund2005}.

Low-order-refined (LOR) preconditioning avoids high-order matrix assembly by preconditioning the spectral-element operator with a lowest-order discretization on its Gauss--Lobatto sub-grid \cite{BelloMaldonadoFischer2019,Canuto1994,CanutoGervasioQuarteroni2010,FrancoCamierAndrejPazner2020,Orszag1980Spectral,Pazner2020LOR}. For the conforming $\Hone$ subspace of $hp$-DG methods, Pazner and Kolev \cite[Section~3.4]{PaznerKolev2022} introduced the variational-restriction LOR operator $P^\top A_{\lorr,\loc} P$. They proved that sharing the high-order conforming prolongation $P$ lifts element-wise spectral equivalence globally across hanging nodes and degree jumps without trace inequalities. Using closed Gauss--Lobatto and open histopolation bases \cite{Gerritsma2011}, Dohrmann \cite{Dohrmann2021} and Pazner, Kolev, and Dohrmann \cite{PaznerKolevDohrmann2023} established element-wise spectral equivalence across $\Hcurl$ and $\Hdiv$. Subsequently, Pazner, Kolev, and Camier \cite{PaznerKolevCamier2023} developed uniform-degree batched LOR assembly with local topological differential matrices and vertex-coordinate evaluation.

Extending this paradigm to general three-dimensional \emph{variable-order} nonconforming de~Rham complexes leaves three concrete gaps. First, on mixed-order hanging faces in three dimensions, pairwise master--slave face interpolation does not prevent internal bisecting edges from carrying lower degrees than the master face. Enforcing conformity without unnecessarily lowering the entire master face therefore requires global trace compatibility elimination across multiple faces. Second, uniform-degree batched LOR setups \cite{PaznerKolevCamier2023} pass nodal coordinates to hypre's auxiliary-space solvers \cite{FalgoutYang2002,KolevVassilevski2009,KolevVassilevski2012ADS}, which reconstruct nodal-to-edge/face interpolants from entrywise absolute values $|G_\true|$ and lowest-order curl incidence. Across mixed-degree or high-order hanging interfaces ($r\ge 2$), negative Gauss--Lobatto Lagrange weights cause constant-field interpolation mismatch that degrades convergence, particularly in $\Hdiv$. Third, uniform-degree batched LOR kernels in MFEM \cite{AndersonEtAl2021MFEM,PaznerKolevCamier2023} are instantiated only for conforming uniform spaces up to $p=5$, leaving heterogeneous sub-grid assembly, owner-entity restrictions, and variable-order vector partial assembly incomplete.

We address these gaps in MFEM \cite{AndersonEtAl2021MFEM} around three main contributions:
\begin{enumerate}
\item \emph{Three-dimensional mixed-order hanging-interface constraint completion.} We combine iterative minimum-rule propagation, parallel internal-edge compatibility elimination ($C P_0\mathbf{u}_{\mathrm{root}}=\mathbf{0}$), owner-entity restriction $R_X$, and nonblocking degree consensus. This preserves a globally conforming discrete subcomplex without over-constraining master faces (\cref{sec:nc}, \cref{alg:constraints}).
\item \emph{Explicit sparse auxiliary-space interpolants and unified variable-order LOR analysis.} Generalizing \cite{PaznerKolev2022,PaznerKolevCamier2023,PaznerKolevDohrmann2023} to arbitrary element degrees $p_K$, we pair element-local $\mathcal{O}(p_K^3)$ direct-sparse LOR assembly with explicit signed-weight N\'ed\'elec and corner-quadrature Raviart--Thomas interpolants (\cref{sec:lor}). We prove warped-bilinear constant-flux exactness (\cref{lem:warped-quad}), separate three levels of constant-field reproduction (\cref{rem:curved-rtpi}), and establish exact discrete commutativity (\cref{thm:derham-commute}) and variational-restriction spectral equivalence (\cref{thm:spectral-equiv,cor:composite-complex}).
\item \emph{Parallel order-grouped partial assembly and solver validation.} We extend degree-grouped sum factorization across $\Hone$, $\Hcurl$, and $\Hdiv$ (\cref{sec:pa}). Parallel benchmarks compare our solvers against uniform-degree batched LOR, coordinate-based auxiliary reconstruction, and assembled preconditioners across random degrees, $10^8$ $\alpha$/$\beta$ contrasts, 2-irregular curved meshes, Fichera corners, and up to 1{,}024 cores and 90 million unknowns (\cref{sec:experiments}).
\end{enumerate}

\section{Parallel variable-order nonconforming de~Rham spaces}
\label{sec:spaces}

\subsection{Gauss--Lobatto and histopolation bases}
\label{sec:bases}
Let $\calT_h$ be a distributed $k$-irregular hexahedral mesh of $\Omega\subset\RR^3$, where each element $K=T_K(\hat K)$ is the image of $\hat K=[0,1]^3$ under a polynomial mapping $T_K$ of fixed degree and carries a polynomial degree $p_K\ge 1$. Following \cite{Gerritsma2011,PaznerKolevDohrmann2023}, for $p=p_K$ let $0=\xi_0<\dots<\xi_p=1$ be the Gauss--Lobatto points on $[0,1]$, $\{\ell_i\}_{i=0}^p$ the closed Lagrange basis of $\mathbb{P}_p$, and $\{\psi_j\}_{j=0}^{p-1}$ the open integrated Gauss--Lobatto--Legendre (histopolation) basis of $\mathbb{P}_{p-1}$ defined by
\begin{equation}
\int_{\xi_m}^{\xi_{m+1}} \psi_j(\xi)\,d\xi = \delta_{mj}, \qquad 0 \le j, m \le p-1.
\label{eq:intgll}
\end{equation}
Unisolvence of \eqref{eq:intgll} on $\mathbb{P}_{p-1}$ together with $\int_{\xi_m}^{\xi_{m+1}}\ell_i'\,d\xi=\delta_{i,m+1}-\delta_{i,m}$ yields the exact differentiation identity
\begin{equation}
\ell_i' = \psi_{i-1} - \psi_i, \qquad 0\le i\le p, \qquad \psi_{-1}\equiv\psi_p\equiv 0 .
\label{eq:dl}
\end{equation}

Tensor products on $\hat K$ define the local reference sequence. The scalar space $W_h$ uses $\mathbb{Q}_{p,p,p}$ of dimension $(p+1)^3$ spanned by $\ell_i\otimes\ell_j\otimes\ell_k$. The N\'ed\'elec space $\mathbf{V}_h$ uses $\mathbb{Q}_{p-1,p,p}\times\mathbb{Q}_{p,p-1,p}\times\mathbb{Q}_{p,p,p-1}$ of dimension $3p(p+1)^2$ with one open factor (e.g., $(\psi_i\otimes\ell_j\otimes\ell_k)\hat{\mathbf{e}}_1$). The Raviart--Thomas space $\mathbf{Q}_h$ uses $\mathbb{Q}_{p,p-1,p-1}\times\mathbb{Q}_{p-1,p,p-1}\times\mathbb{Q}_{p-1,p-1,p}$ of dimension $3p^2(p+1)$ with two open factors (e.g., $(\ell_i\otimes\psi_j\otimes\psi_k)\hat{\mathbf{e}}_1$, giving normal traces in $\mathbb{Q}_{p-1,p-1}$). Physical finite element fields are mapped from $\hat K$ via scalar pullback in $W_h$ and covariant and contravariant Piola transformations in $\mathbf{V}_h$ and $\mathbf{Q}_h$. By \eqref{eq:intgll} and \eqref{eq:dl}, the local unknowns represent nodal values, sub-edge tangential line integrals, and sub-face normal fluxes on the Gauss--Lobatto sub-grid, so $\grad W_h\subset\mathbf{V}_h$ and $\curl\mathbf{V}_h\subset\mathbf{Q}_h$ hold element by element.

\subsection{Mixed-order hanging-interface constraint completion}
\label{sec:nc}
Throughout the paper, for each space $X\in\{\Hone,\ND,\RT\}$ we distinguish three objects: the \emph{local vector} $\mathbf{u}_\loc\in\RR^{N_\loc^X}$ of \emph{local unknowns} attached to active entities (vertices, edges, faces, and interiors) on each process, with master and slave entities carrying separate entries; the \emph{true vector} $\mathbf{u}\in\RR^{N_\true^X}$ of independent globally conforming \emph{true unknowns}; and the resulting \emph{finite element field} ($u_h\in W_h$, $\bfv_h\in\mathbf{V}_h$, or $\bfw_h\in\mathbf{Q}_h$). Global conformity is encoded by a parallel prolongation $P_X:\RR^{N_\true^X}\to\RR^{N_\loc^X}$ and a communication-free block-diagonal restriction $R_X:\RR^{N_\loc^X}\to\RR^{N_\true^X}$ satisfying
\begin{equation}
R_X P_X = I .
\label{eq:RP}
\end{equation}
An element restriction $G_E$ (with block $G_{E,K}$) maps local vectors to lexicographically ordered element vectors. Global edges are oriented by increasing vertex index and quadrilateral faces by their eight dihedral orientation flags \cite{CervenyDobrevKolev2019,KinnewigWickBeuchler2025}, whose signed permutations of Gauss--Lobatto and histopolation unknowns compose multiplicatively across element restrictions and $P_X$.

Extending fixed-order nonconforming prolongations \cite{CervenyDobrevKolev2019} to three-dimensional variable-order spaces requires resolving internal-edge degree conflicts on hanging faces, as summarized in \cref{alg:constraints}. During minimum-rule propagation (lines~1--2), every conforming or master face $F$ adopts the minimum degree $p_F$ of elements directly abutting $F$ or its slave sub-faces. Edge degrees propagate down the edge master--slave directed acyclic graph and from faces to their boundary edges. Crucially, when an internal edge $e_{\mathrm{int}}$ inside a master quadrilateral face $F_m$ acquires a lower degree $p_{e_{\mathrm{int}}} < p_{F_m}$ from an element across a perpendicular face, we do \emph{not} lower $p_{F_m}$. Instead (lines~4--8), after building the pairwise prolongation $P_0:\RR^{N_{\mathrm{root}}^X}\to\RR^{N_\loc^X}$ on candidate roots, we assemble trace compatibility rows $C\mathbf{u}_\loc=\mathbf{0}$ for $X\in\{\Hone,\ND\}$ ($\RT$ has no edge unknowns). For each internal edge $e_{\mathrm{int}}$ of a slave sub-face with $p_{e_{\mathrm{int}}}<p_{F_m}$, we map oriented $e_{\mathrm{int}}$ into the reference square of $F_m$, evaluate the trace transfer $T_{e_{\mathrm{int}}\leftarrow F_m}$ onto the degree-$p_{F_m}$ edge variant, and emit one signed equation per degree-$p_{F_m}$ edge unknown $r$:
\begin{equation}
s_r\, u_{\loc,e_{\mathrm{int}}}^{(p_{F_m})}[r] - \sum_c s_c\,(T_{e_{\mathrm{int}}\leftarrow F_m})_{rc}\, u_{\loc,F_m}[c] = 0,
\label{eq:C-general}
\end{equation}
where $s_r,s_c\in\{\pm1\}$ encode global entity orientations. Because $P_0$ already expresses the degree-$p_{F_m}$ edge variant from its degree-$p_{e_{\mathrm{int}}}$ roots, substituting $\mathbf{u}_\loc = P_0\mathbf{u}_{\mathrm{root}}$ equates the degree-$p_{F_m}$ face trace along $e_{\mathrm{int}}$ to the prolonged degree-$p_{e_{\mathrm{int}}}$ edge trace for arbitrary $p_{e_{\mathrm{int}}}<p_{F_m}$.

\begin{algorithm}[t]
\caption{Parallel 3D mixed-order constraint completion and owner restriction.}
\label{alg:constraints}
\begin{algorithmic}[1]
\Require mesh $\calT_h$, degrees $\{p_K\}$, space $X\in\{\Hone,\ND,\RT\}$, rank tolerance $\tau=128\varepsilon_{\mathrm{mach}}$
\State initialize entity degree bitmasks from local elements; exchange updates via nonblocking consensus \cite{HoeflerSiebertLumsdaine2010NBX}
\State propagate minimum degrees iteratively across master--slave faces/edges and from faces to boundary edges until convergence
\State form pairwise master--slave and multi-variant prolongation $P_0:\RR^{N_{\mathrm{root}}^X}\to\RR^{N_\loc^X}$
\If{$X\in\{\Hone,\ND\}$}
  \State assemble signed internal-edge trace rows $C\mathbf{u}_\loc=\mathbf{0}$ via \eqref{eq:C-general} for every internal edge $e_{\mathrm{int}}$ with $p_{e_{\mathrm{int}}}<p_{F_m}$; form $\tilde C = C P_0$ with row scales $\sigma_r = \sum_l\sum_j |C_{rl}(P_0)_{lj}|$
  \State group coupled candidate roots by union--find and route rows to component owners
  \State reduce each row against normalized echelon rows $U_i$ ($\tilde C_{r,\cdot}\gets \tilde C_{r,\cdot}-\mu U_{i,\cdot}$, $\sigma_r\gets\sigma_r+|\mu||U_{i,k}|$); pick strictly largest $|\tilde C_{r,j^*}|$ (ties broken by smallest global root ID); accept iff $|\tilde C_{r,j^*}|>\tau\sigma_r$, drop $|\tilde C_{rj}|\le\tau\sigma_r$, and normalize by $\tilde C_{r,j^*}$
  \State back-substitute in reverse pivot order ($u_{j_i}=-\sum_{k} U_{i,k} u_k$) to form $Q:\RR^{N_\true^X}\to\RR^{N_{\mathrm{root}}^X}$
\Else
  \State $Q \gets I$
\EndIf
\State set $P_X \gets P_0 Q$ and owner-entity restriction $R_X \gets S R_0$ (dropping unused local variants)
\end{algorithmic}
\end{algorithm}

Mathematically, the globally conforming subcomplex is characterized by the exact algebraic nullspace $\ker(C P_0)$: provided slave geometric mappings match restrictions of their master mappings, tangential gradients and surface curls commute with pullbacks and preserve $\ker(C)$, so $W_h\xrightarrow{\grad}\mathbf{V}_h\xrightarrow{\curl}\mathbf{Q}_h$ is a discrete subcomplex with $\mathrm{range}(P_X)$. In floating-point arithmetic, \cref{alg:constraints} enforces $\tilde C\mathbf{u}_{\mathrm{root}}=\mathbf{0}$ up to the scale-relative numerical rank threshold $\tau\sigma_r$ ($\tau=128\varepsilon_{\mathrm{mach}}$), where $\sigma_r$ bounds accumulated magnitude prior to cancellation. Disjoint interface patches remain decoupled across connected components, so redundant equations (e.g., from adjacent sub-faces) are discarded cleanly without a global dense solve, and deterministic global-ID tie-breaking guarantees partition-independent pivot choices.

To illustrate \cref{alg:constraints}, consider two coarse cubes $[0,2]\times[0,1]^2$ in which $[1,2]\times[0,1]^2$ is refined to level~1 and $[1,1.5]\times[\frac12,1]\times[0,\frac12]$ to level~2, yielding $16$ hexahedra. Assigning $p_{K_*}=1$ to the grandchild $K_*=[1.25,1.5]\times[0.5,0.75]\times[0,0.25]$ (away from $x_1=1$) and $p_K=2$ to the remaining $15$ elements lowers the perpendicular level-1 face $[1,1.5]\times\{\frac12\}\times[0,\frac12]$ and its edge $e_{\mathrm{int}}=\{1\}\times\{\frac12\}\times[0,\frac12]$ to degree~$1$, while the coarse master face $F_m=\{1\}\times[0,1]^2$ retains $p_{F_m}=2$. In $\Hone$, let $(u_{1,0},u_{1,1},u_{1,2})$ be the $Q_2$ candidate roots along $y=\frac12$ on $F_m$ at $z=0,\frac12,1$. Equating the quadratic trace at $z=\frac14$ ($\frac38 u_{1,0}+\frac34 u_{1,1}-\frac18 u_{1,2}$) to the linear trace along $e_{\mathrm{int}}$ ($\frac12 u_{1,0}+\frac12 u_{1,1}$) via \eqref{eq:C-general} gives the representative row of $C P_0$:
\begin{equation}
-\tfrac{1}{8}\,u_{1,0} + \tfrac{1}{4}\,u_{1,1} - \tfrac{1}{8}\,u_{1,2} = 0,
\qquad \sigma_r = \tfrac{9}{4}, \qquad \tau\sigma_r = 288\varepsilon_{\mathrm{mach}} \approx 6.39\times 10^{-14}.
\label{eq:ie-row}
\end{equation}
Pivoting on $|\tilde C_{r,(1,1)}|=\frac14 > \tau\sigma_r$ eliminates the face-center root via $u_{1,1}=\frac12(u_{1,0}+u_{1,2})$ while retaining compatible $Q_2$ face modes on $F_m$ such as $(y-\frac12)z(z-\frac12)$. Across this 16-hexahedron patch, root completion maps $(N_\loc^X, N_{\mathrm{root}}^X)\to N_\true^X$ as $(251,153)\to 152$ in $\Hone$, $(618,417)\to 416$ in $\ND$, and $(474,386)\to 386$ in $\RT$.

Because $P_X$ is rectangular, \eqref{eq:RP} alone does not uniquely determine $R_X$ outside $\mathrm{range}(P_X)$. \Cref{alg:constraints} fixes $R_X = S R_0$ explicitly. Each candidate root $j$ is owned by a unique rank and master entity of global minimum degree $r_j$. Row $j$ of $R_0$ vanishes on slave rows; it selects the degree-$r_j$ local unknown with a signed unit entry when stored on the owner rank, and otherwise interpolates locally ($r_j^\loc\to r_j$) from the lowest locally stored variant $r_j^\loc>r_j$. The row selector $S$ keeps surviving free roots ($S Q = I$, hence $R_X P_X = S R_0 P_0 Q = I$). Dynamic degree consensus uses the nonblocking protocol of Hoefler, Siebert, and Lumsdaine \cite{HoeflerSiebertLumsdaine2010NBX} (synchronous nonblocking sends with a nonblocking barrier), avoiding late in-flight messages and global all-to-all communication.

\section{Order-grouped variable-order partial assembly}
\label{sec:pa}
Partial assembly (PA) represents the high-order bilinear form in factorized form \cite{AndersonEtAl2021MFEM,KronbichlerKormann2012}
\begin{equation}
A_\ho^X = P_X^\top G_E^\top B^\top D\, B\, G_E P_X,
\label{eq:pa-factor}
\end{equation}
where $B$ applies sum-factorized basis or derivative evaluations at quadrature points and the block-diagonal tensor $D$ stores quadrature weights, geometric factors, and coefficients. To reuse uniform-degree tensor-product kernels on variable-order meshes, each process partitions its active elements by polynomial degree,
\begin{equation}
\calT_h^\loc = \bigsqcup_{g=1}^{N_g} \calE_g, \qquad p_K = p_g \;\;\text{for all } K \in \calE_g, \qquad p_1 < \dots < p_{N_g}.
\label{eq:pa-groups}
\end{equation}
A grouped element restriction $G_E$ packs each group $\calE_g$ into a contiguous slice of the element vector. Diffusion/mass ($\Hone$), curl--curl/mass ($\Hcurl$), and div--div/mass ($\Hdiv$) integrators then allocate per-group quadrature tensors $D_g$ ($n_{\mathrm{data}}q_g^3$ scalars per element) and apply batched sum factorization $B_g^\top D_g B_g$ group by group between a single global prolongation $G_E P_X$ and its transpose $P_X^\top G_E^\top$.

\begin{proposition}[Complexity of order-grouped partial assembly]
\label{prop:pa-complexity}
Assume bounded-degree element mappings, $\mathcal{O}(1)$ pointwise coefficient evaluation, and $q_K = p_K + c$ quadrature points per direction with fixed $c$. For each $X\in\{\Hone,\ND,\RT\}$, the per-process element-local setup and storage of \eqref{eq:pa-factor}--\eqref{eq:pa-groups} scale as $\sum_{K\in\calT_h^\loc}\mathcal{O}(p_K^3)$, and each operator application requires $\sum_{K\in\calT_h^\loc}\mathcal{O}(p_K^4) + \mathcal{O}(\mathrm{nnz}(P_X))$ operations. By contrast, a fully assembled matrix incurs $\sum_{K\in\calT_h^\loc}\mathcal{O}(p_K^6)$ element-local storage/application and $\sum_{K\in\calT_h^\loc}\mathcal{O}(p_K^9)$ naive quadrature assembly.
\end{proposition}
\begin{proof}
Each element $K\in\calE_g$ stores $n_{\mathrm{data}}q_K^3$ values ($n_{\mathrm{data}}\le 9$), and sum-factorized Jacobian evaluation over a bounded number of geometric nodes costs $\mathcal{O}(q_K^3)=\mathcal{O}(p_K^3)$ operations. During operator evaluation, $B_g$ and $B_g^\top$ execute three one-dimensional contractions per component in $\mathcal{O}(p_K q_K^3)=\mathcal{O}(p_K^4)$ operations per element, $G_E$ and $G_E^\top$ touch each of the $\mathcal{O}(p_K^3)$ element unknowns once, and $P_X,P_X^\top$ cost $\mathcal{O}(\mathrm{nnz}(P_X))$. By contrast, dense element matrices of size $N_K=\mathcal{O}(p_K^3)$ contain $N_K^2=\mathcal{O}(p_K^6)$ entries and require $\mathcal{O}(q_K^3 N_K^2)=\mathcal{O}(p_K^9)$ operations under direct quadrature assembly.
\end{proof}

\section{Direct-sparse variable-order LOR preconditioning}
\label{sec:lor}

\subsection{Gauss--Lobatto sub-grid and variational restriction}
\label{sec:lor-grid}
On each element $K$ of degree $p_K$, the Gauss--Lobatto points partition $\hat K$ into $p_K^3$ sub-cells $\hat\tau_{\mathbf{o}}=[\xi_{o_1},\xi_{o_1+1}]\times[\xi_{o_2},\xi_{o_2+1}]\times[\xi_{o_3},\xi_{o_3+1}]$, $\mathbf{o}\in\{0,\dots,p_K-1\}^3$, placing the local unknowns of $W_h$, $\mathbf{V}_h$, and $\mathbf{Q}_h$ in one-to-one correspondence with the lowest-order $Q_1$, $\ND_1$, and $\RT_0$ unknowns on the sub-grid \cite{Dohrmann2021,PaznerKolevDohrmann2023}. Following the variational-restriction formulation \cite[Section~3.4]{PaznerKolev2022,PaznerKolevCamier2023}, we assemble the LOR operator on the same local space and couple it through the exact same conforming prolongation $P_X$:
\begin{equation}
A_\ho^X = P_X^\top A_{\ho,\loc}^X P_X, \qquad A_\lorr^X = P_X^\top A_{\lorr,\loc}^X P_X ,
\label{eq:global-ops}
\end{equation}
where $A_{\lorr,\loc}^X = G_E^\top\,\mathrm{blkdiag}_K(A_{\lorr,K}^X)\,G_E$ sums the $p_K^3$ lowest-order sub-cell matrices on each element without constructing a globally refined mesh.

\subsection{Direct-sparse heterogeneous LOR assembly}
\label{sec:lor-mat}
Accumulating $A_{\lorr,K}^X$ in a dense $N_K\times N_K$ buffer would reintroduce $\mathcal{O}(p_K^6)$ element-local memory. \Cref{alg:lor} generalizes the uniform-degree macro-element sparse assembly of \cite{PaznerKolevCamier2023} to arbitrary element-wise degrees $p_K$ without compile-time degree limits. It scatters each sub-cell matrix ($8\times8$ for $Q_1$, $12\times12$ for $\ND_1$, $6\times6$ for $\RT_0$, integrated by $2\times2\times2$ Gauss--Lobatto quadrature) directly into a sparse local matrix $A_\loc$. Caching $T_K$ once at the $(p_K+1)^3$ Gauss--Lobatto nodes reduces corner evaluations by up to $8\times$. Mapping sub-cell quadrature points back to parent reference coordinates ensures consistent evaluation of variable coefficients.

\begin{algorithm}[t]
\caption{Direct-sparse variable-order LOR assembly and auxiliary-space setup.}
\label{alg:lor}
\begin{algorithmic}[1]
\Require space $X_h$, $X\in\{\Hone,\ND,\RT\}$; coefficients $\alpha,\beta$; essential true unknowns $\mathcal{I}_{\mathrm{ess}}$
\State initialize sparse $A_\loc\in\RR^{N_\loc^X\times N_\loc^X}$ and reference element $\hat K_1\in\{Q_1,\ND_1,\RT_0\}$ ($2^3$ GLL rule)
\ForAll{$K \in \calT_h^\loc$}
  \If{$p_K = 1$}
    \State assemble $A_K$ directly into $A_\loc$
  \Else
    \State cache $\bfx_{ijk} = T_K(\xi_i,\xi_j,\xi_k)$ for $0\le i,j,k\le p_K$
    \ForAll{sub-cells $\hat\tau_{\mathbf{o}}$ of $K$}
      \State form trilinear sub-cell map from cached corners; evaluate $\alpha,\beta$ at parent coordinates
      \State scatter signed sub-cell matrix $A_\tau$ directly into $A_\loc$
    \EndFor
  \EndIf
\EndFor
\State form $A_\lorr^X \gets P_X^\top A_\loc P_X$ and eliminate essential rows/columns in $\mathcal{I}_{\mathrm{ess}}$
\State set up BoomerAMG ($\Hone$), AMS with \eqref{eq:lor-ops} ($\ND$), or ADS with \eqref{eq:lor-ops} ($\RT$)
\end{algorithmic}
\end{algorithm}

\subsection{Explicit signed-weight auxiliary-space operators}
\label{sec:lor-aux}
To precondition $A_\lorr^{\ND}$ and $A_\lorr^{\RT}$ via AMS and ADS \cite{HiptmairXu2007,KolevVassilevski2009,KolevVassilevski2012ADS}, hypre requires discrete differential operators together with nodal-to-edge/face interpolants. While hypre's coordinate-vector interface is convenient for lowest-order or uniform-degree batched LOR setups \cite{PaznerKolevCamier2023}, it reconstructs nodal edge interpolation via $(\tilde{\Pi}_d\mathbf{q})_e = \frac{1}{2}(G_\true\bfx_d)_e(|G_\true|\mathbf{q})_e$ \cite{KolevVassilevski2009}. Consider an interior sub-edge $e=(v_0\to v_1)$ meeting an interface node $v_0$ constrained to master roots $\{m_j\}$ with Lagrange weights $\{\alpha_j\}$ ($\sum_j\alpha_j=1$). Row $e$ of $G_\true$ then has $+1$ at $v_1$ and $-\alpha_j$ at $m_j$. For linear ($r=1$) hanging nodes, $\alpha_j\ge 0$ gives $\sum_j|\alpha_j|=1$ and $\tilde\Pi_d$ remains exact. Across mixed-degree or $r\ge 2$ hanging interfaces, however, Gauss--Lobatto Lagrange weights change sign ($\sum_j|\alpha_j|>1$), causing constant-field mismatch. For example, on two conforming cubes $[-1,0]\times[0,1]^2$ ($p_1=2$) and $[0,1]^3$ ($p_2=3$), interpolating the $p_2=3$ node $\xi_1=(5-\sqrt{5})/10$ from the $r=2$ master nodes $0,\frac12,1$ yields $\alpha_0=(1+\sqrt{5})/10$, $\alpha_1=4/5$, and $\alpha_2=(1-\sqrt{5})/10<0$. Taking $|G_\true|$ flips the sign of $\alpha_2$, producing $(\tilde{\Pi}_1\mathbf{1})_e = (1+(\sqrt{5}-1)/10)\Delta x_1 \approx 1.1236\,\Delta x_1$ on the edge and $1.2778\,\Delta x_1$ at the tensor-product face node $(\xi_1,\xi_1)$. In ADS, coordinate-based face interpolation additionally assumes $\pm1$ face--edge incidence in $C_\true$, which ceases to hold on constrained faces.

We instead assemble all four auxiliary operators sparsely on local vectors and restrict them to true unknowns via $P_X$ and the owner restriction $R_X$ of \cref{sec:nc}:
\begin{equation}
\begin{aligned}
G_\true &= R_\ND\, G_\loc\, P_\Hone, &\qquad
\Pi_{d,\true}^{\ND} &= R_\ND\, \Pi_{d,\loc}^{\ND}\, P_\Hone, \\
C_\true &= R_\RT\, C_\loc\, P_\ND, &\qquad
\Pi_{d,\true}^{\RT} &= R_\RT\, \Pi_{d,\loc}^{\RT}\, P_\Hone,
\end{aligned}
\qquad d \in \{1,2,3\}.
\label{eq:lor-ops}
\end{equation}
Because shared local unknowns share the same entity, degree variant, orientation, and geometric trace, the local formulas are single-valued. Note a key distinction in \eqref{eq:lor-ops}: by \cref{thm:derham-commute}, $G_\loc\,\mathrm{range}(P_\Hone)\subseteq\mathrm{range}(P_\ND)$ and $C_\loc\,\mathrm{range}(P_\ND)\subseteq\mathrm{range}(P_\RT)$, so $G_\true$ and $C_\true$ are independent of the choice of left inverse in \eqref{eq:RP}. By contrast, $\Pi^X_{d,\loc}P_\Hone\mathbf{q}$ need not lie in $\mathrm{range}(P_X)$ across mixed-order or curved hanging interfaces. Consequently, $\Pi^X_{d,\true}$ relies specifically on the owner-entity restriction $R_X = S R_0$, which preserves signed weights in $P_\Hone$ (recovering $(\Pi^{\ND}_{1,\true}\mathbf{1})_e=\Delta x_1$ on the $(2,3)$ interface above) and yields globally conforming fields upon re-applying $P_X$.

For a sub-edge $e=(v_0\to v_1)$ with orientation sign $s_e\in\{\pm1\}$,
\begin{equation}
(G_\loc \mathbf{q}_\loc)_e = s_e\,(q_{v_1} - q_{v_0}), \qquad
(\Pi_{d,\loc}^{\ND} \mathbf{q}_\loc)_e = \tfrac12\, s_e\,\big((\bfx_{v_1})_d - (\bfx_{v_0})_d\big)\,(q_{v_0} + q_{v_1}).
\label{eq:G-PiND}
\end{equation}
For a sub-face $f$ normal to $\hat{\mathbf{e}}_c$ with oriented boundary sub-edges $(e_0,e_1)$ along $c_1=c+1$ and $(e_2,e_3)$ along $c_2=c+2$ (modulo $3$) and sign $s_f\in\{\pm1\}$, Stokes' theorem gives
\begin{equation}
(C_\loc \bfv_\loc)_f = s_f\,\big( s_{e_0} v_{e_0} + s_{e_3} v_{e_3} - s_{e_1} v_{e_1} - s_{e_2} v_{e_2} \big).
\label{eq:C}
\end{equation}
On the bilinear sub-face surrogate $\bfx_f(s,t)=(1-s)(1-t)\bfx_{00}+s(1-t)\bfx_{10}+st\bfx_{11}+(1-s)t\bfx_{01}$ with corner normals $\bfn_f^{(k)}=\partial_s\bfx_f\times\partial_t\bfx_f$ at $(\hat s_k,\hat t_k)\in\{0,1\}^2$, the $2\times2$ Gauss--Lobatto flux rule yields
\begin{equation}
(\Pi_{d,\loc}^{\RT} \mathbf{q}_\loc)_f = \frac{s_f}{4} \sum_{k=0}^{3} \big(\bfn_f^{(k)}\big)_d \, q_{v_k}.
\label{eq:rt-pi}
\end{equation}
\Cref{tab:derham-summary} summarizes all components across the complex.

\begin{table}[t]
\centering
\caption{Variable-order PA and direct-sparse LOR components ($p=p_K$; local row nonzeros in parentheses).}
\label{tab:derham-summary}
\scriptsize
\setlength{\tabcolsep}{3.5pt}
\begin{tabular}{@{}l l l l l l@{}}
\toprule
Space & Basis & PA forms & Sub-cell & Local aux.\ ops (nnz/row) & Solver \\
\midrule
$\Hone$ & $Q_p$ (GLL) & diff., mass & $Q_1$ ($8{\times}8$) & --- & AMG \\
$\Hcurl$ & $\ND_p$ (GLL/IGLL) & curl--curl, mass & $\ND_1$ ($12{\times}12$) & $G_\loc$, $\Pi_{d,\loc}^{\ND}$ ($\le 2$) & AMS \\
$\Hdiv$ & $\RT_{p-1}$ (GLL/IGLL) & div--div, mass & $\RT_0$ ($6{\times}6$) & $C_\loc$, $\Pi_{d,\loc}^{\RT}$ ($\le 4$); $G_\loc$, $\Pi_{d,\loc}^{\ND}$ & ADS \\
\bottomrule
\end{tabular}
\end{table}

\subsection{Theoretical analysis}
\label{sec:theory}
We now analyze sparsity, constant-field reproduction, discrete commutativity, and spectral equivalence; $\preceq$ denotes the Loewner order.

\begin{proposition}[Complexity and sparsity of direct-sparse LOR]
\label{prop:lor-complexity}
Under the hypotheses of \cref{prop:pa-complexity}, \cref{alg:lor} assembles $A_\loc$ and $G_\loc,C_\loc,\Pi_{d,\loc}^{\ND},\Pi_{d,\loc}^{\RT}$ in $\sum_{K\in\calT_h^\loc}\mathcal{O}(p_K^3)$ operations and memory. The row nonzeros of $A_\loc$ are bounded independently of $\{p_K\}$ by the mesh vertex/edge valence: at most $11$ for $\Hdiv$ on any hexahedral mesh, and on Cartesian-like patches at most $27$ for $\Hone$ and $33$ for $\Hcurl$. Forming $A_\lorr^X = P_X^\top A_\loc P_X$ and \eqref{eq:lor-ops} adds work and interface fill governed by $\mathrm{nnz}(P_X)$.
\end{proposition}
\begin{proof}
Caching $T_K$ at $(p_K+1)^3$ nodes and scattering $p_K^3$ fixed-size sub-cell matrices costs $\mathcal{O}(p_K^3)$ per element without dense element buffers. Prior to applying $P_X$, master and slave entities carry distinct local unknowns. Each sub-face belongs to at most two sub-cells ($6+6-1=11$ faces in $\Hdiv$). On interior and Cartesian patches, sub-vertices belong to at most $8$ sub-cells ($27$ vertices in $\Hone$) and sub-edges to at most $4$ sub-cells ($33$ edges in $\Hcurl$), while general boundary entities are bounded by mesh valence. By \eqref{eq:G-PiND}--\eqref{eq:rt-pi}, local auxiliary rows have at most $2$ or $4$ nonzeros; $p$-dependent interface fill enters only through $P_X$ and $R_X$.
\end{proof}

\begin{lemma}[Constant-flux reproduction on warped bilinear sub-faces]
\label{lem:warped-quad}
For any bilinear sub-face $f$ parametrized by $\bfx_f(s,t)$ with corners $\bfx_{00},\bfx_{10},\bfx_{11},\bfx_{01}$ and corner normals $\bfn_f^{(k)}$,
\begin{equation}
\frac{1}{4}\sum_{k=0}^{3} \bfn_f^{(k)}
= \int_{[0,1]^2} \partial_s \bfx_f \times \partial_t \bfx_f\,ds\,dt
= \frac{1}{2}(\bfx_{11} - \bfx_{00}) \times (\bfx_{01} - \bfx_{10}).
\label{eq:warped-quad-id}
\end{equation}
Hence, for every constant vector $\mathbf{c}\in\RR^3$ and nodal vector $\mathbf{1}$, $\sum_{d=1}^{3} c_d (\Pi_{d,\loc}^{\RT} \mathbf{1}_\loc)_f = s_f \int_f \mathbf{c}\cdot \bfn\,dS$, even when the four corners of $f$ are non-coplanar.
\end{lemma}
\begin{proof}
Write $\bfx_f(s,t) = \bfx_{00} + s\mathbf{a} + t\mathbf{b} + st\bfw$ with $\mathbf{a}=\bfx_{10}-\bfx_{00}$, $\mathbf{b}=\bfx_{01}-\bfx_{00}$, and $\bfw=\bfx_{11}-\bfx_{10}-\bfx_{01}+\bfx_{00}$. Since $\partial_s \bfx_f = \mathbf{a} + t\bfw$, $\partial_t \bfx_f = \mathbf{b} + s\bfw$, and $\bfw\times\bfw = \mathbf{0}$,
\[
\partial_s \bfx_f \times \partial_t \bfx_f
= \mathbf{a}\times\mathbf{b} + s(\mathbf{a}\times\bfw) + t(\bfw\times\mathbf{b})
\]
is affine in $(s,t)$, so the $2\times 2$ trapezoidal rule is exact and integrates it to $\mathbf{a}\times\mathbf{b} + \frac12(\mathbf{a}\times\bfw + \bfw\times\mathbf{b})$. Expanding $\frac12(\bfx_{11}-\bfx_{00})\times(\bfx_{01}-\bfx_{10})=\frac12(\mathbf{a}+\mathbf{b}+\bfw)\times(\mathbf{b}-\mathbf{a})$ yields the identical vector, proving \eqref{eq:warped-quad-id}; contracting with $s_f\mathbf{c}$ gives the flux identity.
\end{proof}

\begin{remark}[Constant-field reproduction across affine and curved meshes]
\label{rem:curved-rtpi}
For a constant vector $\mathbf{c}\in\RR^3$ and nodal vector $\mathbf{1}$ ($P_\Hone\mathbf{1}=\mathbf{1}_\loc$), let $\mathbf{d}^X_\loc(\mathbf{c})$ denote the exact local degrees of freedom of $\mathbf{c}$ and $\mathcal{I}_X$ the map from true vectors to finite element fields. Three distinct properties hold:
\begin{enumerate}
\item[\textup{(P1)}] \emph{Local functional exactness:} $\sum_d c_d(\Pi^{\ND}_{d,\loc}\mathbf{1}_\loc)_e = \mathbf{d}^{\ND}_\loc(\mathbf{c})_e$ holds on every smooth curved sub-edge since $\int_e\mathbf{c}\cdot d\bfx=\mathbf{c}\cdot(\bfx_{v_1}-\bfx_{v_0})$, while $\sum_d c_d(\Pi^{\RT}_{d,\loc}\mathbf{1}_\loc)_f$ is exact on warped bilinear sub-faces (\cref{lem:warped-quad}) and approximates higher-degree curved sub-faces.
\item[\textup{(P2)}] \emph{True-unknown projection agreement:} applying $R_\ND$ yields $\sum_d c_d\Pi^{\ND}_{d,\true}\mathbf{1}=R_\ND\mathbf{d}^{\ND}_\loc(\mathbf{c})$ identically on every mesh and degree distribution.
\item[\textup{(P3)}] \emph{Conforming prolongation and physical constant reproduction:} re-applying $P_\ND$ gives $P_\ND\sum_d c_d\Pi^{\ND}_{d,\true}\mathbf{1}=\mathbf{d}^{\ND}_\loc(\mathbf{c})$ and $\mathcal{I}_\ND(\sum_d c_d\Pi^{\ND}_{d,\true}\mathbf{1})\equiv\mathbf{c}$ if and only if $\mathbf{c}\in\mathbf{V}_h$; in $\RT$, $\mathcal{I}_\RT(\sum_d c_d\Pi^{\RT}_{d,\true}\mathbf{1})\equiv\mathbf{c}$ holds if and only if $\mathbf{c}\in\mathbf{Q}_h$ and $R_\RT(\sum_d c_d\Pi^{\RT}_{d,\loc}\mathbf{1}_\loc-\mathbf{d}^{\RT}_\loc(\mathbf{c}))=\mathbf{0}$.
\end{enumerate}
On affine hexahedra, all three properties hold globally for every $p_K\ge 1$. On curved $Q_q$ ($q\ge 2$) hexahedra, covariant Piola pullbacks $J_K^\top\mathbf{c}$ have degrees up to $(q-1,q,q)$, so $\min_K p_K\ge q$ suffices to ensure $\mathbf{c}\in\mathbf{V}_h$ for all $\mathbf{c}\in\RR^3$ \cite{ArnoldBoffiBonizzoni2015}. For example, under the quadratic shear $T(s,t,r)=(s,t+\gamma s(1-s),r)$ on $[0,1]^3$ with $\mathbf{c}=(0,1,0)^\top$ and $p_K=1$, $\Pi^{\ND}_{2,\true}\mathbf{1}$ matches all edge integrals, yet $\mathcal{I}_\ND(\Pi^{\ND}_{2,\true}\mathbf{1})=(-\gamma(1-2s),1,0)^\top\neq\mathbf{c}$ ($L^2$ error $|\gamma|/\sqrt{3}$). Across a curved $r=1$ hanging face, $P_\ND\Pi^{\ND}_{2,\true}\mathbf{1}$ enforces tangential continuity rather than matching unconstrained half-edge integrals ($\pm\gamma/4$). Taking $p_K\ge 2$ restores $\mathcal{I}_\ND(\Pi^{\ND}_{2,\true}\mathbf{1})\equiv\mathbf{c}$, whereas $G_\true$ and $C_\true$ remain exact for all $p_K\ge 1$ on curved meshes (\cref{thm:derham-commute}).
\end{remark}

\begin{theorem}[Exact discrete de~Rham commutativity]
\label{thm:derham-commute}
Let $W_h$, $\mathbf{V}_h$, and $\mathbf{Q}_h$ be the variable-order spaces of \cref{sec:spaces} on a possibly curved, $k$-irregular hexahedral mesh with common degree distribution $\{p_K\}$, compatible interface geometric traces, and prolongations/restrictions satisfying \eqref{eq:RP} such that $W_h \xrightarrow{\grad} \mathbf{V}_h \xrightarrow{\curl} \mathbf{Q}_h$ forms a globally conforming discrete subcomplex characterized by $\mathrm{range}(P_X)$. Let $\mathcal{I}_\ND,\mathcal{I}_\RT$ map true vectors to finite element fields. Then \eqref{eq:lor-ops} satisfies:
\begin{enumerate}
  \item $P_\ND\, G_\true = G_\loc\, P_\Hone$ and $P_\RT\, C_\true = C_\loc\, P_\ND$;
  \item $C_\true G_\true = 0$, $\mathcal{I}_\ND G_\true\mathbf{q}=\grad u_h$, and $\mathcal{I}_\RT C_\true\mathbf{v}=\curl\bfv_h$ for the fields $u_h\in W_h$, $\bfv_h\in\mathbf{V}_h$ associated with $\mathbf{q},\mathbf{v}$, so $\curl(\mathcal{I}_\ND G_\true\mathbf{q})=\mathbf{0}$ and $\divg(\mathcal{I}_\RT C_\true\mathbf{v})=0$ pointwise on every element;
  \item for unconstrained curl--curl and div--div stiffness matrices $S^{\ND},S^{\RT}$ assembled via \eqref{eq:global-ops} with $\alpha>0$ and positive quadrature weights,
  \[
  \mathrm{range}(G_\true) \subseteq \ker(S_{\ho}^{\ND}) \cap \ker(S_{\lorr}^{\ND}), \quad
  \mathrm{range}(C_\true) \subseteq \ker(S_{\ho}^{\RT}) \cap \ker(S_{\lorr}^{\RT}),
  \]
  and the same inclusions hold for the boundary-eliminated operators on compatible homogeneous essential boundary subspaces.
\end{enumerate}
\end{theorem}
\begin{proof}
\emph{Claim 1.} For $\mathbf{q}\in\RR^{N_\true^\Hone}$ representing $u_h\in W_h$ with local vector $P_\Hone\mathbf{q}$, the pullbacks $\hat u = u_h \circ T_K$ and $\widehat{\grad u_h} = J^\top ((\grad u_h)\circ T_K)$ ($J=DT_K$) satisfy $\widehat{\grad u_h}=\hat\grad\hat u$. By \eqref{eq:intgll}--\eqref{eq:dl}, the sub-edge functional of $\grad u_h$ on $e=(v_0\to v_1)$ is $s_e(\hat u(v_1)-\hat u(v_0))=(G_\loc P_\Hone\mathbf{q})_e$, so $G_\loc P_\Hone\mathbf{q}$ is the local vector of $\grad u_h\in\mathbf{V}_h$. By the subcomplex property, $G_\loc P_\Hone\mathbf{q}=P_\ND\mathbf{g}$ for some $\mathbf{g}$; applying $R_\ND$ and \eqref{eq:RP} gives $\mathbf{g}=G_\true\mathbf{q}$ and $P_\ND G_\true\mathbf{q}=G_\loc P_\Hone\mathbf{q}$. Analogously, for $\bfv_h\in\mathbf{V}_h$ with local vector $P_\ND\mathbf{v}$, the Piola identity $\widehat{\curl\bfv_h}=\hat\curl\hat\bfv$ and Stokes' theorem show that $(C_\loc P_\ND\mathbf{v})_f$ is the local vector of $\curl\bfv_h\in\mathbf{Q}_h$, so $C_\loc P_\ND\mathbf{v}\in\mathrm{range}(P_\RT)$ and left-multiplying by $R_\RT$ yields $P_\RT C_\true\mathbf{v}=C_\loc P_\ND\mathbf{v}$.

\emph{Claim 2.} Every sub-face corner bounds two edges entering \eqref{eq:C} with opposite signs, giving $C_\loc G_\loc=0$ and thus $C_\true G_\true = R_\RT C_\loc P_\ND G_\true = R_\RT C_\loc G_\loc P_\Hone = 0$. The field identities follow directly from Claim~1 and $\curl\grad=0$, $\divg\curl=0$.

\emph{Claim 3.} With positive weights $\{\omega_q\}$ and $\alpha>0$, $\mathbf{w}^\top S_{\ho,K}^\ND \mathbf{w} = \sum_q \omega_q \alpha(\bfx_q)|\curl \bfw_h(\bfx_q)|^2 \ge 0$, so $\mathbf{w}^\top S_{\ho,K}^\ND\mathbf{w}=0\implies S_{\ho,K}^\ND\mathbf{w}=\mathbf{0}$. Since $\mathcal{I}_\ND G_\true\mathbf{q}=\grad u_h$ is curl-free, $S_{\ho,\loc}^\ND P_\ND G_\true\mathbf{q}=\mathbf{0}$ and $S_\ho^\ND G_\true\mathbf{q}=\mathbf{0}$. On each sub-cell $\tau$, \eqref{eq:G-PiND} restricts $G_\loc P_\Hone\mathbf{q}$ to the $\ND_1$ edge functionals of the trilinear nodal interpolant's gradient, which is curl-free on $\tau$, giving $S_\lorr^\ND G_\true\mathbf{q}=P_\ND^\top S_{\lorr,\loc}^\ND G_\loc P_\Hone\mathbf{q}=\mathbf{0}$. Identical reasoning using $\divg\curl=0$ and $\ND_1\xrightarrow{\curl}\RT_0$ proves the $C_\true$ inclusions, and homogeneous essential boundary values vanish on identity-replaced rows and columns.
\end{proof}

\begin{remark}
\Cref{thm:derham-commute} ensures that AMS and ADS receive exact discrete gradient and curl images ($\mathrm{range}(G_\true)\subseteq\ker(S^\ND)$, $\mathrm{range}(C_\true)\subseteq\ker(S^\RT)$) on arbitrary variable-order nonconforming meshes; bounded convergence as $\beta\to 0$ also depends on subspace stability and algebraic multigrid coarsening and smoothing \cite{BakerFalgoutKolevYang2011,HiptmairXu2007}.
\end{remark}

\begin{assumption}[Uniform element-wise spectral equivalence]
\label{ass:local}
For each $X\in\{\Hone,\ND,\RT\}$ and $K\in\calT_h$, let $S_{\ho,K}^X,M_{\ho,K}^X$ and $S_{\lorr,K}^X,M_{\lorr,K}^X$ denote the unit-coefficient high-order and direct-sparse LOR element stiffness and mass matrices. There exist constants $0<\underline{s}\le\overline{s}$ and $0<\underline{m}\le\overline{m}$, independent of $h_K$ and $p_K$, such that
\[
\underline{s}\, S_{\lorr,K}^X \preceq S_{\ho,K}^X \preceq \overline{s}\, S_{\lorr,K}^X, \qquad
\underline{m}\, M_{\lorr,K}^X \preceq M_{\ho,K}^X \preceq \overline{m}\, M_{\lorr,K}^X .
\]
\end{assumption}

On affine parallelepipeds, \cref{ass:local} holds with constants depending only on shape regularity: for $\Hone$ by \cite{Canuto1994,CanutoGervasioQuarteroni2010} and for $\Hcurl,\Hdiv$ in the GLL/IGLL bases by \cite{Dohrmann2021,PaznerKolevDohrmann2023}. Furthermore, whenever the stiffness matrix $S_{\bullet,K}^X$ and the next-space mass matrix $M_{\bullet,K}^{X'}$ ($\bullet\in\{\ho,\lorr\}$) are assembled with the same element or sub-cell quadrature rule---or integrated exactly on affine parallelepipeds---the discrete differential inclusions $\grad W_K\subset\mathbf{V}_K$, $\curl\mathbf{V}_K\subset\mathbf{Q}_K$, and $\divg\mathbf{Q}_K\subset L^2_K$ yield the factorizations
\begin{equation}
S_{\bullet,K}^{\Hone} = D_0^\top M_{\bullet,K}^{\ND} D_0, \qquad
S_{\bullet,K}^{\ND} = D_1^\top M_{\bullet,K}^{\RT} D_1, \qquad
S_{\bullet,K}^{\RT} = D_2^\top M_{\bullet,K}^{L^2} D_2,
\label{eq:stiff-inc}
\end{equation}
through the shared sub-grid incidence matrices $D_0,D_1,D_2$. Under this quadrature-compatibility condition, Loewner congruence transfers mass-matrix equivalence directly to stiffness matrices; under mismatched quadrature rules, stiffness equivalence does not follow from \eqref{eq:stiff-inc} and must be verified directly. On curved elements, where \cref{alg:lor} uses trilinear sub-cell surrogates and $2^3$ GLL quadrature, \cref{ass:local} requires uniform map non-degeneracy, scaled metric equivalence between $DT_K$ and sub-cell Jacobians, and quadrature norm equivalence. Combining \cref{ass:local} with the variational-restriction lifting of \cite[Section~3.4]{PaznerKolev2022} yields global spectral equivalence across the complex.

\begin{theorem}[Uniform spectral equivalence of variable-order nonconforming LOR]
\label{thm:spectral-equiv}
Let $\calT_h$ be a $k$-irregular hexahedral mesh with arbitrary degrees $p_K\ge 1$ and element-wise constant coefficients $\alpha|_K=\alpha_K>0$, $\beta|_K=\beta_K>0$, and let \cref{ass:local} hold. For $X\in\{\Hone,\ND,\RT\}$, let $A_\ho^X$ and $A_\lorr^X$ be assembled via \eqref{eq:global-ops} with identical identity-row/column elimination of essential true unknowns, and set $c_0=\min\{\underline{s},\underline{m},1\}$, $c_1=\max\{\overline{s},\overline{m},1\}$. Then
\begin{equation}
c_0 \, \mathbf{x}^\top A_\lorr^X \mathbf{x}
\le \mathbf{x}^\top A_\ho^X \mathbf{x}
\le c_1 \, \mathbf{x}^\top A_\lorr^X \mathbf{x}
\qquad \text{for all } \mathbf{x} \in \RR^{N_\true^X},
\label{eq:spectral-equiv}
\end{equation}
so $\kappa((A_\lorr^X)^{-1}A_\ho^X)\le c_1/c_0$, with $c_0,c_1$ independent of $h_K$, irregularity $k$, degree distribution and neighbor degree jumps, $(\alpha_K,\beta_K)$, and the parallel partition.
\end{theorem}
\begin{proof}
Scaling the stiffness and mass Loewner bounds in \cref{ass:local} by the positive weights $\alpha_K>0$ and $\beta_K>0$ and adding them yields
\begin{equation}
c_0\, A_{\lorr,K}^X \preceq A_{\ho,K}^X \preceq c_1\, A_{\lorr,K}^X ,
\label{eq:local-equiv}
\end{equation}
uniformly in $K,h_K,p_K,\alpha_K,\beta_K$. Splitting any $\mathbf{x}=\mathbf{x}_0+\mathbf{x}_{\mathrm{ess}}$ into free and essential true unknowns gives $\mathbf{x}^\top A^X\mathbf{x}=\mathbf{x}_0^\top\tilde A^X\mathbf{x}_0+\|\mathbf{x}_{\mathrm{ess}}\|_2^2$ for $A^X\in\{A_\ho^X,A_\lorr^X\}$, where $\tilde A^X$ is the unconstrained operator \eqref{eq:global-ops}. Setting $\mathbf{v}_K=G_{E,K}P_X\mathbf{x}_0$, shared prolongation and element restriction yield
\[
\mathbf{x}_0^\top \tilde A_\ho^X \mathbf{x}_0 = \sum_{K\in\calT_h} \mathbf{v}_K^\top A_{\ho,K}^X \mathbf{v}_K,
\qquad
\mathbf{x}_0^\top \tilde A_\lorr^X \mathbf{x}_0 = \sum_{K\in\calT_h} \mathbf{v}_K^\top A_{\lorr,K}^X \mathbf{v}_K .
\]
Summing \eqref{eq:local-equiv} over $K\in\calT_h$ and adding $\|\mathbf{x}_{\mathrm{ess}}\|_2^2$ with $c_0\le 1\le c_1$ proves \eqref{eq:spectral-equiv} without trace inequalities or degree grading; Courant--Fischer gives $\kappa((A_\lorr^X)^{-1}A_\ho^X)\le c_1/c_0$.
\end{proof}

\begin{remark}[Anisotropic box elements]
\label{rem:aspect}
On axis-aligned boxes, each component mass block in $\Hone,\ND,\RT,L^2$ factors into a shared geometric scalar times a Kronecker product of three 1D SPD mass matrices, making \cref{ass:local}---and hence \eqref{eq:stiff-inc}---independent of element aspect ratio; algebraic multigrid sensitivity on stretched sub-grids is tested in \cref{sec:num-s1b}.
\end{remark}

\begin{corollary}[Composite auxiliary-space preconditioner]
\label{cor:composite-complex}
Under the hypotheses of \cref{thm:spectral-equiv}, let $B^X$ be an SPD preconditioner satisfying $\gamma_0 (A_\lorr^X)^{-1} \preceq B^X \preceq \gamma_1 (A_\lorr^X)^{-1}$ with $0<\gamma_0\le\gamma_1$. Then $\mathrm{spec}(B^X A_\ho^X)\subset[\gamma_0 c_0,\gamma_1 c_1]$ and $\kappa(B^X A_\ho^X)\le (\gamma_1/\gamma_0)(c_1/c_0)$; in exact arithmetic, PCG reduces the $A_\ho^X$-norm of the error by a factor $\varepsilon\in(0,1)$ within $\lceil\tfrac12\sqrt{\kappa(B^X A_\ho^X)}\,\ln(2/\varepsilon)\rceil$ iterations.
\end{corollary}
\begin{proof}
Inverting \eqref{eq:spectral-equiv} gives $c_1^{-1}(A_\lorr^X)^{-1}\preceq(A_\ho^X)^{-1}\preceq c_0^{-1}(A_\lorr^X)^{-1}$ and thus $\gamma_0 c_0(A_\ho^X)^{-1}\preceq B^X\preceq\gamma_1 c_1(A_\ho^X)^{-1}$. Similarity of $B^X A_\ho^X$ to $(A_\ho^X)^{1/2}B^X(A_\ho^X)^{1/2}$ bounds its Rayleigh quotients in $[\gamma_0 c_0,\gamma_1 c_1]$, and standard CG theory \cite{HestenesStiefel1952,Saad2003} yields the iteration bound.
\end{proof}

\begin{remark}
While auxiliary-space theory guarantees mesh-independent $\gamma_0,\gamma_1$ on shape-regular lowest-order grids \cite{HiptmairXu2007,KolevVassilevski2009,KolevVassilevski2012ADS}, Gauss--Lobatto sub-cells exhibit $\mathcal{O}(p)$ boundary aspect ratios and practical solvers apply a single V-cycle, so $p$-robustness of the composite solver is assessed numerically for $p\le 6$ below \cite{Pazner2020LOR,PaznerKolevDohrmann2023}.
\end{remark}

\section{Numerical experiments}
\label{sec:experiments}

\subsection{Experimental setup and timing conventions}
\label{sec:num-setup}
All benchmarks were run on dual-socket AMD EPYC 9755 nodes (256 cores and 768~GB memory per node) with GCC~14.2, OpenMPI~4.1.8 (one rank per core), MFEM \cite{AndersonEtAl2021MFEM}, and hypre~3.1.0 \cite{FalgoutYang2002}. Linear systems are solved by PCG \cite{HestenesStiefel1952,Saad2003} with relative tolerance $10^{-10}$, enforcing true Euclidean relative residual $\|A\mathbf{x}-\mathbf{b}\|_2/\allowbreak\|\mathbf{b}\|_2 \le 10^{-10}$ (except where noted in \cref{tab:s3}), using HMIS coarsening (strength threshold $0.5$, one aggressive coarsening level) and $\ell_1$-Gauss--Seidel smoothing in BoomerAMG, AMS, and ADS \cite{BakerFalgoutKolevYang2011}.

We compare five solver configurations: PA+LOR (order-grouped PA with \cref{alg:lor} and explicit $\Pi^{\ND},\allowbreak\Pi^{\RT}$), PA+LOR with coordinate reconstruction (Coord.\ $\Pi$, replacing explicit $\Pi^{\ND},\allowbreak\Pi^{\RT}$ by hypre's coordinate interface on identical $A_\lorr^X,G_\true,C_\true$), upstream uniform-degree PA+\allowbreak BatchedLOR \cite{PaznerKolevCamier2023} ($p\le 5$, with PA+\allowbreak LegacyLOR at $p=6$), FA+LOR, and assembled FA+Aux (FA+AMG, FA+AMS, or FA+ADS). We separately time finite element space, constraint, and $P_X,\allowbreak R_X$ construction ($t_{\mathrm{fes}}$). Total solver time is reported as $t_{\mathrm{tot}} = t_{\mathrm{op}} + t_{\mathrm{pc}} + t_{\mathrm{sol}}$, combining operator setup, preconditioner setup (including initial hypre setup), and PCG solve; $t_{\mathrm{fes}}$ is excluded because space setup amortizes across solves, and is reported alongside $t_{\mathrm{tot}}$.

\begin{figure}[t]
\centering
\includegraphics[width=\textwidth]{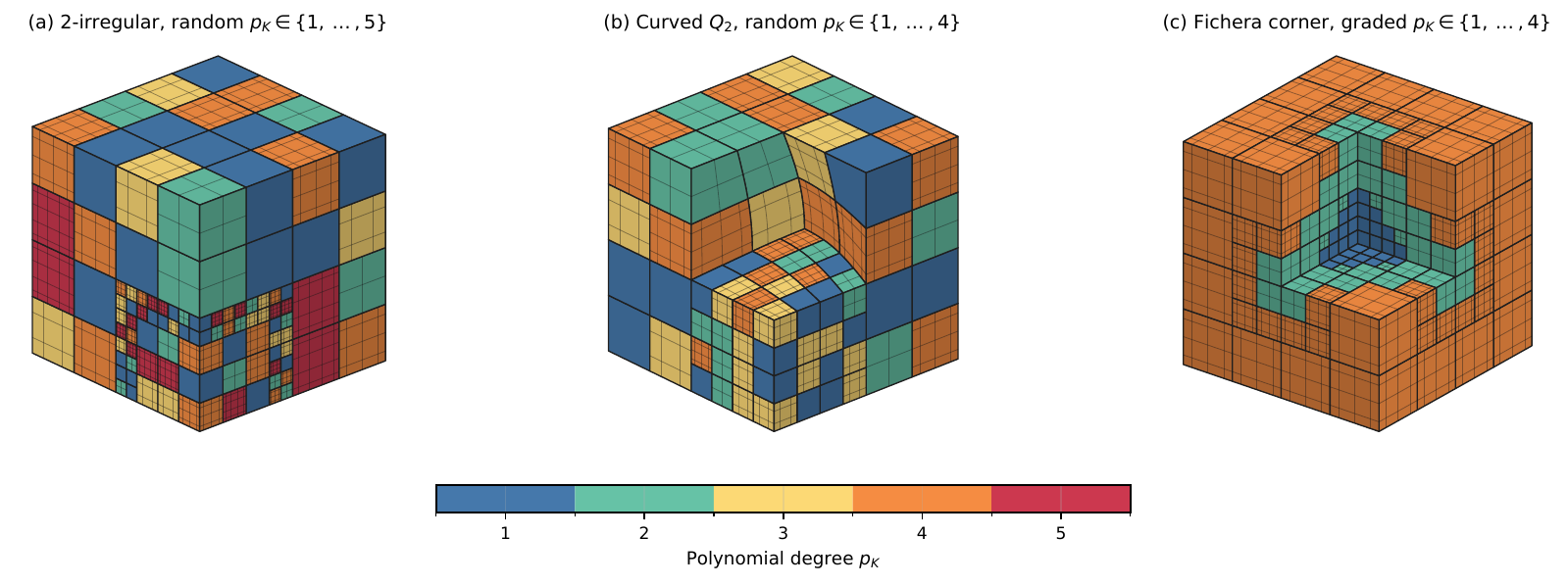}
\caption{Representative 3D variable-order nonconforming meshes colored by element degree $p_K$ with Gauss--Lobatto LOR sub-grids: (a) 2-irregular mesh, random $p_K\in\{1,\dots,5\}$; (b) cutaway curved $Q_2$ ($\delta=0.08$) 1-irregular mesh, random $p_K\in\{1,\dots,4\}$; (c) Fichera corner, 2-level nonconforming graded $p_K\in\{1,\dots,4\}$.}
\label{fig:meshes}
\end{figure}

On $\Omega=[0,1]^3$, 1-irregular meshes refine $[0,0.5]^3$ with 2:1 face balance, and 2-irregular meshes refine the adjacent fine layer without balance to create $16{:}1$ hanging faces. Curved quadratic ($Q_2$) meshes nodally interpolate
\begin{equation}
\mathbf{y}(\bfx) = \bfx + \delta \sin(\pi x_1)\sin(\pi x_2)\sin(\pi x_3)\,(1,\,-0.5,\,0.7)^\top
\label{eq:curved-warp}
\end{equation}
on the coarse grid prior to refinement, with sampled normalized Jacobians $\det(DT_K)/\allowbreak V_{K,0}$ ($9^3$ closed GLL points per element, rounded outward) in $[0.80,1.20]$ ($\delta=0.06$) and $[0.74,1.26]$ ($\delta=0.08$). Smooth accuracy tests use the manufactured fields $u=\sin(\pi x_1)\allowbreak\sin(\pi x_2)\allowbreak\sin(\pi x_3)$ in $\Hone$, $\bfu=(\sin\pi x_2\sin\pi x_3,\allowbreak\,\sin\pi x_3\sin\pi x_1,\allowbreak\,\sin\pi x_1\sin\pi x_2)^\top$ in $\Hcurl$, and $\bfu=\grad u/\pi$ in $\Hdiv$. Representative variable-order meshes are shown in \cref{fig:meshes}.

\subsection{Verification of constraints, commutativity, and interpolants}
\label{sec:num-verify}
To verify the discrete theory of \cref{sec:spaces,sec:lor}, we test $4$- and $64$-rank runs across 1-irregular graded and 2-irregular random meshes ($p_K\in\{1,\dots,4\}$ and $\{2,\dots,4\}$, affine and curved $\delta=0.06$), as well as the minimal 16-hexahedron patch of \cref{sec:nc}. Across all cases, $G_\true$ and $C_\true$ match reference operators, satisfy $P_\ND G_\true=G_\loc P_\Hone$, $P_\RT C_\true=C_\loc P_\ND$, and $C_\true G_\true=0$, reproduce curved curls in $L^2$, yield conforming interface traces for $P_X\Pi^X_{d,\true}\mathbf{q}$, and match PA against FA (with shared curved quadrature $q_K=p_{\max}+3$) to within $1.7\times 10^{-14}$. On the 16-hexahedron internal-edge patch (tested serially for pre-compaction counts $(251,153)\to 152$, $(618,417)\to 416$, $(474,386)\to 386$ and partitioned across MPI ranks), \cref{alg:constraints} yields matching global $N_\true^X$, direct $\Hone,\ND,\RT$ interface jumps below $1.2\times 10^{-13}$, and $Q_2$ master-face mode reproduction to $6.4\times 10^{-17}$. Constant-field tests confirm \cref{rem:curved-rtpi} (including the exact $p_K=1$ error $|\gamma|/\sqrt{3}$ and half-edge difference $\gamma/4$, restored to rounding error for $\min_K p_K\ge 2$), the $27.8\%$ $|G_\true|$ mismatch on $(p_1,p_2)=(2,3)$, and curved $Q_2$ $\Pi^{\RT}$ convergence ($N=4,8,16$ at $p=2$, $\delta=0.06$) with face-flux rates $4.03,4.01$ and $L^2$ rates $1.99,2.00$.

\subsection{Uniform degree scaling and comparison with batched LOR}
\label{sec:num-s1a}
\begin{figure}[t]
\centering
\includegraphics[width=\textwidth]{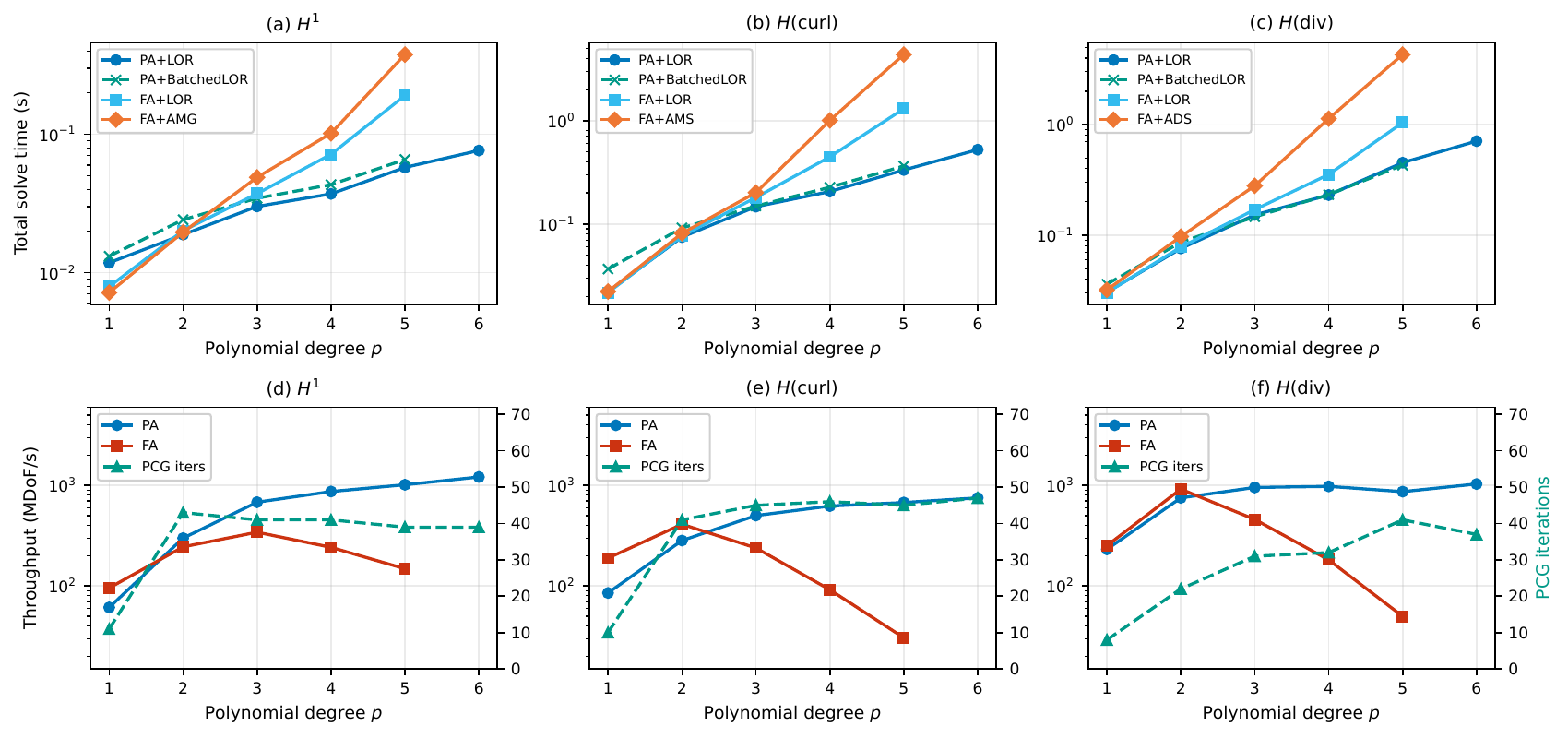}
\caption{Uniform degree scaling ($p=1,\dots,6$) on an $8^3$ mesh (64 cores): (a--c) total solve time $t_{\mathrm{tot}}$; (d--f) operator throughput (solid) and PA+LOR iterations (dashed).}
\label{fig:s1a}
\end{figure}

We first evaluate $p$-robustness and baseline efficiency on a uniform $8^3$ mesh ($p=1,\dots,6$, $\alpha=\beta=1$, 64 cores; \cref{tab:s1a} and \cref{fig:s1a}). Discretization $L^2$ errors agree across solvers to at least four digits for $p\le 5$ (and within $2.5\%$ at $p=6$ in $\Hdiv$ near solver tolerance). PCG iteration counts for PA+LOR remain essentially flat from $p=2$ to $p=6$ across all three spaces, whereas FA+Aux requires fewer iterations on the exact high-order matrix but incurs steep $\mathcal{O}(p^6)$ setup and memory growth.

Comparing PA+LOR against upstream PA+BatchedLOR \cite{PaznerKolevCamier2023} confirms that our runtime heterogeneous assembly matches the iteration counts and CPU solve times of order-templated uniform kernels for $p=2,\dots,5$ (using fewer iterations at $p=1$, where $2^3$ GLL quadrature matches PA). At $p=6$, where BatchedLOR is uninstantiated and upstream MFEM falls back to mesh-refining LegacyLOR, PA+LOR extends seamlessly with lower solve time. Space and prolongation setup $t_{\mathrm{fes}}$ remains a small fraction of $t_{\mathrm{tot}}$ across all degrees, and PA+LOR achieves up to $13\times$ end-to-end speedup over FA+Aux at $p=5$ (full $p=1,\dots,6$ metrics appear in \cref{tab:sm-order}).

\begin{table}[t]
\centering
\caption{Uniform degree scaling ($p=4,5,6$) on an $8^3$ hexahedral mesh (64 cores). Space/prolongation setup $t_{\mathrm{fes}}$, preconditioner setup $t_{\mathrm{pc}}$, and total solve time $t_{\mathrm{tot}}=t_{\mathrm{op}}+t_{\mathrm{pc}}+t_{\mathrm{sol}}$ are in seconds; throughput (MV) is in $\text{MDoF/s}$. Full $p=1,\dots,6$ data (including NNZ and memory) appear in Table~SM1.}
\label{tab:s1a}
\scriptsize
\setlength{\tabcolsep}{3.0pt}
\begin{tabular}{@{}l r r l r r r r r@{}}
\toprule
Space & $p$ & DoFs & Method & It & $t_{\mathrm{fes}}$ & $t_{\mathrm{pc}}$ & $t_{\mathrm{tot}}$ & MV \\
\midrule
$\Hone$ & 4 & 35{,}937 & \textbf{PA+LOR-AMG} & 41 & 0.005 & 0.017 & \textbf{0.037} & \textbf{869} \\
 &  &  & PA+BatchedLOR & 41 & 0.002 & 0.018 & 0.043 & 698 \\
 &  &  & FA+AMG & 21 & 0.005 & 0.031 & 0.102 & 249 \\
 & 5 & 68{,}921 & \textbf{PA+LOR-AMG} & 39 & 0.006 & 0.027 & \textbf{0.058} & \textbf{1{,}016} \\
 &  &  & PA+BatchedLOR & 40 & 0.002 & 0.029 & 0.066 & 943 \\
 &  &  & FA+LOR-AMG & 39 & 0.006 & 0.027 & 0.191 & 149 \\
 &  &  & FA+AMG & 30 & 0.006 & 0.115 & 0.377 & 152 \\
 & 6 & 117{,}649 & \textbf{PA+LOR-AMG} & 39 & 0.010 & 0.037 & \textbf{0.076} & \textbf{1{,}213} \\
 &  &  & PA+LegacyLOR & 40 & 0.002 & 0.054 & 0.097 & 1{,}150 \\
\midrule
$\Hcurl$ & 4 & 104{,}544 & \textbf{PA+LOR-AMS} & 46 & 0.019 & 0.053 & \textbf{0.206} & \textbf{624} \\
 &  &  & PA+BatchedLOR & 47 & 0.002 & 0.058 & 0.227 & 594 \\
 &  &  & FA+AMS & 20 & 0.018 & 0.190 & 1.003 & 93 \\
 & 5 & 201{,}720 & \textbf{PA+LOR-AMS} & 45 & 0.041 & 0.086 & \textbf{0.333} & \textbf{676} \\
 &  &  & PA+BatchedLOR & 46 & 0.003 & 0.090 & 0.364 & 690 \\
 &  &  & FA+LOR-AMS & 45 & 0.041 & 0.084 & 1.300 & 31 \\
 &  &  & FA+AMS & 27 & 0.041 & 0.567 & 4.344 & 34 \\
 & 6 & 345{,}744 & \textbf{PA+LOR-AMS} & 47 & 0.085 & 0.128 & \textbf{0.523} & \textbf{752} \\
 &  &  & PA+LegacyLOR & 47 & 0.003 & 0.167 & 0.589 & 732 \\
\midrule
$\Hdiv$ & 4 & 101{,}376 & \textbf{PA+LOR-ADS} & 32 & 0.006 & 0.085 & \textbf{0.231} & \textbf{979} \\
 &  &  & PA+BatchedLOR & 33 & 0.002 & 0.073 & 0.232 & 992 \\
 &  &  & FA+ADS & 12 & 0.006 & 0.416 & 1.131 & 190 \\
 & 5 & 196{,}800 & \textbf{PA+LOR-ADS} & 41 & 0.007 & 0.151 & \textbf{0.453} & \textbf{867} \\
 &  &  & PA+BatchedLOR & 39 & 0.002 & 0.120 & 0.430 & 887 \\
 &  &  & FA+LOR-ADS & 41 & 0.007 & 0.136 & 1.041 & 50 \\
 &  &  & FA+ADS & 15 & 0.008 & 1.343 & 4.298 & 50 \\
 & 6 & 338{,}688 & \textbf{PA+LOR-ADS} & 37 & 0.010 & 0.259 & \textbf{0.709} & \textbf{1{,}033} \\
 &  &  & PA+LegacyLOR & 41 & 0.002 & 0.211 & 0.743 & 739 \\
\bottomrule
\end{tabular}
\end{table}

\subsection{Random variable-order distributions and auxiliary interpolation}
\label{sec:num-random}
To test ungraded degree jumps and isolate the solver impact of signed-weight auxiliary interpolation, we assign independent random degrees $p_K$ drawn uniformly from $\{1,\dots,p_{\max}\}$ across conforming, 1- and 2-irregular, and curved $Q_2$ meshes (\cref{tab:random} and \cref{fig:s1b}(d)). With explicit $\Pi^{\ND},\Pi^{\RT}$, PCG iterations remain tightly bounded across abrupt degree jumps, multiple random seeds, and curved interfaces. Iterations rise only moderately on direct 2-irregular ($16{:}1$ sub-face) interfaces in both PA+LOR and FA+Aux, consistent with AMG coarsening across wider master--slave couplings; PA+LOR achieves up to $4\times$ speedup over FA+Aux at $p_{\max}=5$.

Comparing PA+LOR with explicit $(\Pi^{\ND},\Pi^{\RT})$ against coordinate-based reconstruction (Coord.\ $\Pi$ on the identical $A_\lorr^X,G_\true,C_\true$) reveals a clear structural contrast between $\Hcurl$ and $\Hdiv$. In AMS ($\Hcurl$), exact $G_\true$ already captures discrete gradients (\cref{thm:derham-commute}) and point relaxation damps localized edge mismatch, so coordinate reconstruction adds at most $0\text{--}3$ iterations across all tested regimes. In ADS ($\Hdiv$), however, replacing both explicit $\Pi^{\ND}$ and $\Pi^{\RT}$ with coordinate reconstruction builds edge and face interpolants from $|G_\true|$ and non-incidence constrained rows of $C_\true$. This inflates PCG iterations by $18\%\text{--}64\%$ across constant-coefficient random mixed-order meshes ($12\%$ under $10^6$ $(\alpha,\beta)$-jumps; $33\to 54$ and $+31\%$ in $t_{\mathrm{tot}}$ on conforming $p_K\in\{1,\dots,5\}$) and by $75\%$ on the graded Fichera mesh ($36\to 63$ in \cref{sec:num-s2}). Together, the explicit $(\Pi^{\ND},\Pi^{\RT})$ interpolants substantially improve $\Hdiv$ solver efficiency across mixed-order interfaces.

\begin{table}[t]
\centering
\caption{Random variable-order benchmarks ($p_K \in \{1,\dots,p_{\max}\}$, 64 cores except bottom block on 512 cores) comparing PA+LOR with explicit signed-weight interpolation (Expl.\ $\Pi$), coordinate-based reconstruction (Coord.\ $\Pi$), and assembled FA+Aux. Multi-seed, low-$\beta$, and high-contrast comparisons appear in Tables~SM2--SM4.}
\label{tab:random}
\scriptsize
\setlength{\tabcolsep}{2.8pt}
\begin{tabular}{@{}l l r r r r r r r@{}}
\toprule
Mesh ($p_K$ range) & Sp. & DoFs & \multicolumn{3}{c}{PCG Iters} & \multicolumn{3}{c}{$t_{\mathrm{tot}}$ (s)} \\
\cmidrule(lr){4-6} \cmidrule(lr){7-9}
 &  &  & Expl.$\Pi$ & Coord.$\Pi$ & FA+Aux & Expl.$\Pi$ & Coord.$\Pi$ & FA+Aux \\
\midrule
Conf.\ $8^3$ ($1\dots 5$) & $\Hone$ & 18{,}565 & \textbf{32} & --- & 25 & \textbf{0.038} & --- & 0.095 \\
 & $\Hcurl$ & 60{,}015 & \textbf{42} & 42 & 24 & \textbf{0.262} & 0.272 & 1.014 \\
 & $\Hdiv$ & 64{,}389 & \textbf{33} & 54 & 13 & \textbf{0.350} & 0.458 & 1.227 \\
\midrule
1-irreg.\ ($1\dots 5$) & $\Hone$ & 12{,}878 & \textbf{34} & --- & 29 & \textbf{0.036} & --- & 0.084 \\
 & $\Hcurl$ & 42{,}440 & \textbf{47} & 47 & 29 & \textbf{0.237} & 0.237 & 0.943 \\
 & $\Hdiv$ & 46{,}367 & \textbf{39} & 52 & 15 & \textbf{0.338} & 0.367 & 1.183 \\
\midrule
2-irreg.\ ($1\dots 4$) & $\Hone$ & 18{,}108 & \textbf{39} & --- & 32 & \textbf{0.038} & --- & 0.069 \\
 & $\Hcurl$ & 62{,}953 & \textbf{54} & 55 & 32 & \textbf{0.279} & 0.286 & 0.738 \\
 & $\Hdiv$ & 72{,}121 & \textbf{41} & 53 & 18 & \textbf{0.368} & 0.386 & 0.850 \\
\midrule
Curved 1-irreg.\ ($1\dots 4$) & $\Hone$ & 8{,}148 & \textbf{30} & --- & 23 & \textbf{0.029} & --- & 0.043 \\
 & $\Hcurl$ & 26{,}675 & \textbf{41} & 42 & 24 & \textbf{0.178} & 0.185 & 0.338 \\
 & $\Hdiv$ & 28{,}967 & \textbf{33} & 42 & 13 & \textbf{0.235} & 0.243 & 0.506 \\
\midrule
512c 1-irreg.\ ($1\dots 4$) & $\Hone$ & 16.82M & \textbf{49} & --- & --- & \textbf{10.79} & --- & --- \\
 & $\Hcurl$ & 29.36M & \textbf{55} & --- & --- & \textbf{29.78} & --- & --- \\
 & $\Hdiv$ & 33.84M & \textbf{43} & --- & --- & \textbf{49.32} & --- & --- \\
\bottomrule
\end{tabular}
\end{table}

\subsection{Robustness to vanishing reaction, contrast, anisotropy, and curvature}
\label{sec:num-s1b}
\begin{figure}[t]
\centering
\includegraphics[width=0.98\textwidth]{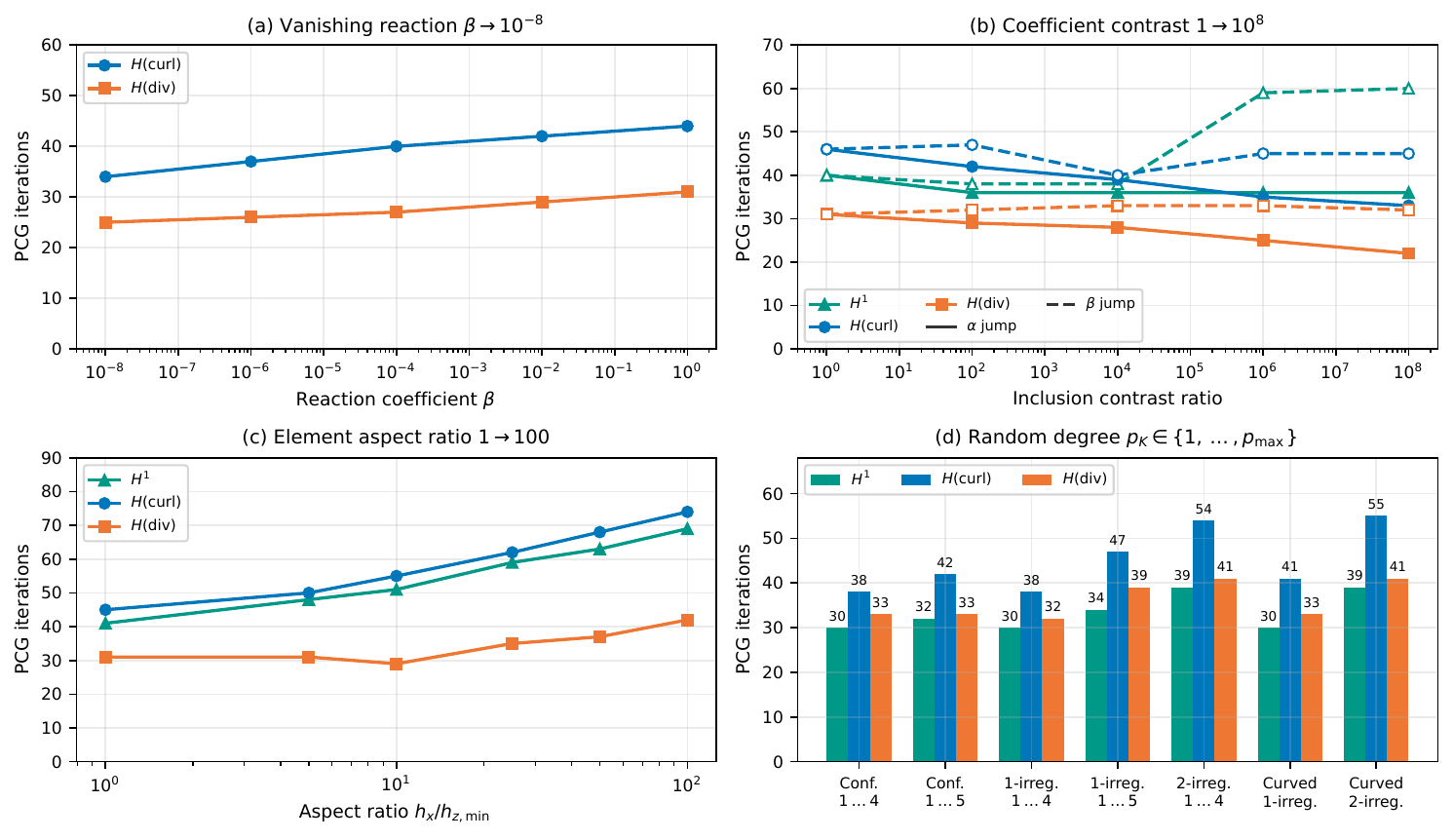}
\caption{PA+LOR robustness on variable-order nonconforming meshes (64 cores): (a) vanishing reaction $\beta\to 10^{-8}$; (b) $[0.25,0.75]^3$ contrast up to $10^8$ in $\alpha$ (solid) and $\beta$ (dashed); (c) element aspect ratio up to $100$; (d) random degrees $p_K\in\{1,\dots,p_{\max}\}$.}
\label{fig:s1b}
\end{figure}

\Cref{fig:s1b}(a--c) evaluates three parameter and geometric regimes on 64 cores. First, as $\beta\to 10^{-8}$ (with algebraically consistent $\mathbf{b}=A\mathbf{u}_{\mathrm{proj}}$), LOR-AMS and LOR-ADS iterations remain bounded and decrease slightly, consistent with the discrete kernel inclusions of \cref{thm:derham-commute}. Second, for an element-aligned inclusion $[0.25,0.75]^3$ with contrast $C\in[1,10^8]$ applied separately to $\alpha$, to $\beta$, or simultaneously, positive-weight scaling and element-wise summation preserve the uniform spectral equivalence bounds $(c_0,c_1)$ of \cref{thm:spectral-equiv}. Iterations remain flat across $\alpha$-jumps in all spaces and across $\beta$-jumps in $\Hcurl,\Hdiv$, showing only a bounded increase in $\Hone$ when the consistent $Q_1$ mass block dominates inside the inclusion. Third, geometric $x_3$-stretching up to aspect ratio $100$ induces only mild logarithmic iteration growth consistent with scalar AMG coarsening on stretched sub-grids (\cref{rem:aspect}). Finally, shared-quadrature curved $Q_2$ ($\delta=0.08$) 1-irregular runs match FA discretization errors and achieve $1.5\text{--}2.2\times$ speedups over FA+Aux on random $p_K\in\{1,\dots,4\}$ (\cref{tab:random}) and $2.3\text{--}4.4\times$ on graded $p_K\in\{1,\dots,4\}$.

\subsection{Fichera corner with graded \texorpdfstring{$hp$}{hp} meshes}
\label{sec:num-s2}
To assess combined geometric and polynomial grading near singularities \cite{Schwab1998}, we discretize $[-1,1]^3\setminus(0,1]^3$ by refining seven unit cubes $R\in\{2,3,4\}$ times (Levels~2--4), adding two 2:1-balanced nonconforming corner refinements ($|\bfx|<1.05/2^{\ell-1}$), and grading $p_K\in\{1,2,4\}$ radially (\cref{fig:meshes}(c) and \cref{tab:s2}). On Level~3 (64 cores), PA+LOR achieves $3.6\text{--}4.1\times$ speedup over FA+Aux across all three spaces while reducing peak memory by up to $2.9\times$. Moreover, the combined explicit $(\Pi^{\ND},\Pi^{\RT})$ auxiliary interpolants cut $\Hdiv$ iterations from $63$ (Coord.\ $\Pi$, $12.8\text{\,s}$) to $36$ ($8.8\text{\,s}$). Comparing FA+LOR against FA+Aux and PA+LOR shows that sparse LOR preconditioning and matrix-free operator evaluation contribute complementary savings, enabling Level~4 solves approaching $15$ million unknowns on 64 cores with nearly unchanged iteration counts.

\begin{table}[t]
\centering
\caption{Fichera corner benchmark (Levels~3--4, 64 cores) with 2-level nonconforming refinement and graded $p_K \in \{1,2,4\}$. Times are in seconds ($t_{\mathrm{fes}}=0.54\text{\,s}$ in $\Hcurl$ and $0.09\text{\,s}$ in $\Hdiv$ on Level~3) and aggregate peak RSS in GiB. Level~2 appears in Table~SM5.}
\label{tab:s2}
\scriptsize
\setlength{\tabcolsep}{2.8pt}
\begin{tabular}{@{}l l r l r r r r r@{}}
\toprule
Sp. & Lvl & DoFs & Method & It & NNZ (M) & $t_{\mathrm{pc}}$ & $t_{\mathrm{tot}}$ & RSS \\
\midrule
$\Hone$ & L3 & 577{,}458 & \textbf{PA+LOR-AMG} & 45 & 13.4 (LOR) & 0.28 & \textbf{0.64} & \textbf{6.7} \\
 &  &  & FA+LOR-AMG & 45 & 97.3 & 0.26 & 1.89 & 9.4 \\
 &  &  & FA+AMG & 23 & 97.3 & 0.26 & 2.62 & 9.3 \\
 & L4 & 4.81M & \textbf{PA+LOR-AMG} & 47 & 94.3 (LOR) & 3.83 & \textbf{6.68} & \textbf{12.5} \\
\midrule
$\Hcurl$ & L3 & 1.77M & \textbf{PA+LOR-AMS} & 53 & 69.1 (LOR) & 1.84 & \textbf{7.76} & \textbf{10.0} \\
 &  &  & PA+LOR (Coord.$\Pi$) & 54 & 69.1 (LOR) & 1.91 & 7.91 & 10.0 \\
 &  &  & FA+LOR-AMS & 53 & 617.6 & 1.69 & 14.19 & 29.3 \\
 &  &  & FA+AMS & 21 & 617.6 & 4.14 & 27.64 & 29.2 \\
 & L4 & 14.57M & \textbf{PA+LOR-AMS} & 56 & 495.3 (LOR) & 26.36 & \textbf{67.41} & \textbf{31.0} \\
\midrule
$\Hdiv$ & L3 & 1.80M & \textbf{PA+LOR-ADS} & 36 & 27.1 (LOR) & 2.93 & \textbf{8.81} & \textbf{11.8} \\
 &  &  & PA+LOR (Coord.$\Pi$) & 63 & 27.1 (LOR) & 2.96 & 12.84 & 11.6 \\
 &  &  & FA+LOR-ADS & 36 & 425.8 & 2.70 & 12.27 & 21.7 \\
 &  &  & FA+ADS & 15 & 425.8 & 10.74 & 32.95 & 28.6 \\
 & L4 & 14.73M & \textbf{PA+LOR-ADS} & 37 & 192.7 (LOR) & 43.94 & \textbf{84.07} & \textbf{41.9} \\
\bottomrule
\end{tabular}
\end{table}

\subsection{Parallel scaling}
\label{sec:num-s3}
\begin{figure}[t]
\centering
\includegraphics[width=\textwidth]{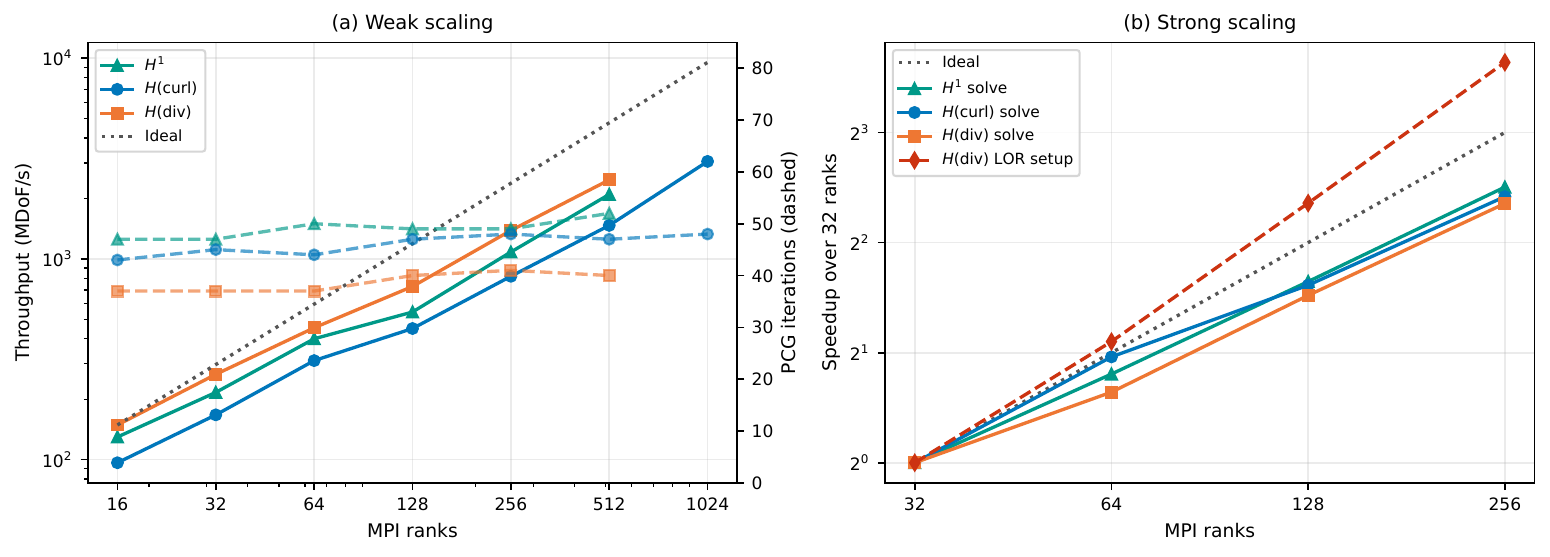}
\caption{Parallel scaling of PA+LOR on 1-irregular graded $hp$ meshes: (a) weak scaling up to 1{,}024 cores; (b) single-node strong scaling from 32 to 256 cores.}
\label{fig:s3}
\end{figure}

\Cref{tab:s3} and \cref{fig:s3} evaluate weak and strong scaling of PA+LOR on 1-irregular graded meshes ($p_K\in\{2,3,4\}$) across 256-core nodes. Under weak scaling up to 512 cores for $\Hone$ ($28.7\text{M}$ unknowns) and $\Hdiv$ ($42.6\text{M}$ unknowns) and 1{,}024 cores for $\Hcurl$ ($90.6\text{M}$ unknowns), PCG iterations and per-rank memory remain nearly flat while matrix-free operator throughput scales nearly linearly past $3\text{\,BDoF/s}$. Operator setup $t_{\mathrm{op}}$ remains small relative to $t_{\mathrm{sol}}$, whereas $t_{\mathrm{sol}}$ grows moderately across multiple nodes due to inter-node communication in the AMG, AMS, and ADS V-cycles. In single-node strong scaling from 32 to 256 cores ($8\times$ cores), $t_{\mathrm{sol}}$ achieves $5.1\text{--}5.7\times$ speedup, while purely local sub-grid assembly $t_{\mathrm{lor}}$ scales superlinearly as per-rank sub-grid working sets shrink.

\begin{table}[t]
\centering
\caption{Parallel weak scaling ($16\to 512$ and $1{,}024$ cores) and single-node strong scaling ($32\to 256$ cores) of PA+LOR on 1-irregular graded meshes ($p_K\in\{2,3,4\}$). Times are in seconds, MV in $\text{MDoF/s}$, and per-rank memory ($\mathrm{RSS}_{\max}$) in MiB. The $\Hone$ and $\Hdiv$ series enforce true relative residual $\le 10^{-10}$; $\Hcurl$ uses preconditioned tolerance $10^{-10}$ (true residual $\le 8\times 10^{-9}$). Full rank counts appear in Table~SM6.}
\label{tab:s3}
\scriptsize
\setlength{\tabcolsep}{2.1pt}
\begin{tabular}{@{}l l r r r r r r r r@{}}
\toprule
Sp. & Mode & Ranks & DoFs & It & $t_{\mathrm{op}}$ & $t_{\mathrm{lor}}$ & $t_{\mathrm{sol}}$ & MV & RSS$_{\max}$ \\
\midrule
$\Hone$ & Weak & 16 & 0.91M & 47 & 0.11 & 1.06 & 1.33 & 130 & 215 \\
 & Weak & 64 & 3.67M & 50 & 0.13 & 1.60 & 2.35 & 400 & 249 \\
 & Weak & 256 (1N) & 14.24M & 49 & 0.19 & 2.39 & 3.03 & 1{,}081 & 267 \\
 & Weak & \textbf{512 (2N)} & \textbf{28.67M} & \textbf{52} & \textbf{0.19} & \textbf{3.35} & \textbf{8.13} & \textbf{2{,}093} & \textbf{271} \\
 & Strong & 32$\to$256 & 7.23M & 49$\to$49 & 0.48$\to$0.09 & 10.81$\to$1.04 & 7.93$\to$1.40 & 197$\to$1{,}435 & 610$\to$196 \\
\midrule
$\Hcurl$ & Weak & 16 & 1.20M & 43 & 0.09 & 1.19 & 7.40 & 96 & 311 \\
 & Weak & 64 & 5.82M & 44 & 0.15 & 2.96 & 12.85 & 311 & 422 \\
 & Weak & 256 (1N) & 22.73M & 48 & 0.25 & 5.17 & 22.67 & 818 & 487 \\
 & Weak & 512 (2N) & 44.95M & 47 & 0.27 & 5.51 & 46.63 & 1{,}467 & 494 \\
 & Weak & \textbf{1{,}024 (4N)} & \textbf{90.58M} & \textbf{48} & \textbf{0.25} & \textbf{6.04} & \textbf{87.15} & \textbf{3{,}056} & \textbf{507} \\
 & Strong & 32$\to$256 & 5.82M & 44$\to$46 & 0.24$\to$0.12 & 5.47$\to$0.63 & 24.70$\to$4.62 & 176$\to$859 & 591$\to$230 \\
\midrule
$\Hdiv$ & Weak & 16 & 1.31M & 37 & 0.11 & 0.79 & 9.67 & 149 & 454 \\
 & Weak & 64 & 5.57M & 37 & 0.14 & 1.63 & 19.14 & 453 & 542 \\
 & Weak & 256 (1N) & 21.62M & 41 & 0.18 & 2.17 & 29.68 & 1{,}384 & 568 \\
 & Weak & \textbf{512 (2N)} & \textbf{42.62M} & \textbf{40} & \textbf{0.25} & \textbf{2.30} & \textbf{41.27} & \textbf{2{,}472} & \textbf{581} \\
 & Strong & 32$\to$256 & 5.57M & 36$\to$38 & 0.27$\to$0.06 & 3.48$\to$0.28 & 29.81$\to$5.83 & 246$\to$1{,}571 & 916$\to$257 \\
\bottomrule
\end{tabular}
\end{table}

\section{Conclusions}
\label{sec:conclusions}
We have developed a parallel matrix-free solver framework for three-dimensional variable-order nonconforming $\Hone$, $\Hcurl$, and $\Hdiv$ discretizations that unites mixed-order hanging-interface constraint completion, explicit signed-weight auxiliary-space interpolants on heterogeneous Gauss--Lobatto sub-grids, and order-grouped partial assembly. Sharing the completed conforming prolongations across the high-order, LOR, and auxiliary operators preserves exact discrete de~Rham commutativity and global variational-restriction spectral equivalence without degree grading.

Parallel experiments up to $p=6$, 1{,}024 cores, and 90 million unknowns show that the framework matches templated uniform-degree batched LOR for $p=2,\dots,5$ on conforming grids while remaining robust across ungraded random degree jumps, 2-irregular curved interfaces, vanishing reaction terms, and $10^8$ $\alpha$/$\beta$ contrasts. In particular, preserving signed constraint weights in the explicit $(\Pi^{\ND},\Pi^{\RT})$ auxiliary interpolants substantially reduces $\Hdiv$ iteration counts compared with coordinate-based reconstruction on mixed-order meshes. Mild iteration growth on high-aspect-ratio boundary layers is consistent with algebraic multigrid coarsening rather than LOR equivalence (\cref{rem:aspect}).

Natural next steps include GPU acceleration of the degree-grouped PA and heterogeneous sparse LOR kernels \cite{KolevCEED2021,PaznerKolevCamier2023} and integration with goal-oriented $hp$-adaptive refinement loops \cite{MelenkWohlmuth2001}.

\paragraph{Code and data availability}
The variable-order constraint, partial assembly, and LOR implementations have been developed in MFEM \cite{AndersonEtAl2021MFEM} for upstream contribution. A self-contained reproducibility archive containing the MFEM source patch, standalone benchmark drivers, archived CSV datasets, and figure/table scripts accompanies this article as arXiv ancillary files (\texttt{anc/}) alongside the Supplementary Material appended below.

\FloatBarrier

\clearpage
\appendix
\setcounter{section}{0}
\setcounter{table}{0}
\setcounter{figure}{0}
\renewcommand{\thesection}{SM\arabic{section}}
\renewcommand{\thetable}{SM\arabic{table}}
\renewcommand{\thefigure}{SM\arabic{figure}}
\renewcommand{\theHsection}{SM\arabic{section}}
\renewcommand{\theHtable}{SM\arabic{table}}
\renewcommand{\theHfigure}{SM\arabic{figure}}

\begin{center}
  {\large\bfseries\MakeUppercase{Supplementary Material}}
\end{center}
\vspace{0.5em}

\section{Reproducibility package and timing conventions}
\label{sec:sm-repro}
All code patches, standalone benchmark drivers, archived CSV datasets, and plotting scripts are bundled inside a single self-contained archive (\texttt{anc/}) with no external repository dependencies:
\begin{enumerate}
  \item \texttt{patches/} and \texttt{data/}: a unified MFEM source patch against public commit \texttt{b370a2c36ef8} (yielding \texttt{a8b2a9c266c5}, hypre~3.1.0) and all archived CSV/diagnostic datasets;
  \item \texttt{src/} and \texttt{run\_local\_checks.sh}: standalone C++17 drivers and a four-rank MPI suite verifying discrete commutativity ($C_{\true}G_{\true}=0$), curved $Q_2$ curl reproduction, owner-restricted interpolation traces, constant-field reproduction, the serial and MPI-partitioned 16-hexahedron patch of \cref{sec:nc}, and solver smoke tests across $\Hone$, $\Hcurl$, and $\Hdiv$;
  \item \texttt{scripts/} and \texttt{reproduce\_paper\_artifacts.sh}: an automated pipeline regenerating all main-text and supplementary LaTeX tables and vector figures directly from \texttt{data/}.
\end{enumerate}

Across all tables, $t_{\mathrm{fes}}$ denotes space initialization, constraint completion, and $P_X,R_X$ construction; $t_{\mathrm{op}}$ denotes PA/FA operator setup; $t_{\mathrm{lor}}$ denotes sparse LOR assembly; $t_{\mathrm{pc}}$ denotes total preconditioner setup (including $t_{\mathrm{lor}}$, auxiliary operators, and initial hypre setup); $t_{\mathrm{sol}}$ denotes PCG iterations; and $t_{\mathrm{tot}} = t_{\mathrm{op}} + t_{\mathrm{pc}} + t_{\mathrm{sol}}$ denotes solver wall time ($t_{\mathrm{fes}}$ reported separately). Throughput (MV) is in $\text{MDoF/s}$, aggregate RSS in $\text{GiB}$, and per-rank $\mathrm{RSS}_{\max}$ in $\text{MiB}$.

\section{Extended benchmark tables}
\label{sec:sm-tables}
\Cref{tab:sm-order,tab:sm-pi,tab:sm-seeds,tab:sm-contrast,tab:sm-fichera-l2,tab:sm-scaling} provide the complete parameter sweeps summarized in \cref{sec:experiments}:
\begin{itemize}
  \item \cref{tab:sm-order} reports the complete uniform degree sweep $p=1,\dots,6$ on an $8^3$ mesh (64 cores), comparing PA+LOR against upstream \texttt{BatchedLORAssembly} ($p\le 5$), \texttt{LegacyAssembleSystem} ($p=6$), FA+LOR, and assembled FA+Aux;
  \item \cref{tab:sm-pi} compares explicit signed-weight auxiliary interpolation ($\Pi^{\ND},\Pi^{\RT}$) against coordinate-based reconstruction (Coord.\ $\Pi$) on the identical LOR operator across 12 conforming, nonconforming, curved, vanishing-reaction, high-contrast, and Fichera regimes;
  \item \cref{tab:sm-seeds} reports independent random degree realizations ($s\in\{0,1,2\}$) with exact element-per-degree histograms and curved $Q_2$ 2-irregular runs;
  \item \cref{tab:sm-contrast} reports PCG iterations under separate principal ($\alpha$), reaction ($\beta$), and simultaneous ($\alpha=\beta$) inclusion jumps up to $10^8$;
  \item \cref{tab:sm-fichera-l2} reports the coarse Level~2 Fichera corner results complementing Levels~3--4 in the main text;
  \item \cref{tab:sm-scaling} reports all tested weak-scaling ($16\to 512$ and $1{,}024$ ranks) and strong-scaling ($32\to 256$ ranks) points.
\end{itemize}

\begin{table}[p]
\centering
\caption{Complete uniform degree sweep ($p=1,\dots,6$, $8^3$ mesh, 64 cores) across $\Hone$, $\Hcurl$, and $\Hdiv$ (RSS in GiB, MV in MDoF/s).}
\label{tab:sm-order}
\tiny
\renewcommand{\arraystretch}{0.88}
\setlength{\tabcolsep}{3.0pt}
\begin{tabular}{@{}l r r l r r r r r r@{}}
\toprule
Sp. & $p$ & DoFs & Method & It & NNZ (M) & $t_{\mathrm{fes}}$ & $t_{\mathrm{tot}}$ & MV & RSS \\
\midrule
$\Hone$ & 1 & 729 & PA+LOR-AMG & 11 & 0.02 & 0.005 & 0.012 & 61 & 4.53 \\
 &  &  & PA+BatchedLOR & 31 & 0.02 & 0.002 & 0.013 & 70 & 4.63 \\
 &  &  & FA+LOR-AMG & 11 & 0.02 & 0.004 & 0.008 & 96 & 4.52 \\
 &  &  & FA+AMG & 11 & 0.02 & 0.004 & 0.007 & 93 & 4.52 \\
 & 2 & 4{,}913 & PA+LOR-AMG & 43 & 0.03 & 0.004 & 0.019 & 299 & 4.92 \\
 &  &  & PA+BatchedLOR & 43 & 0.12 & 0.002 & 0.024 & 252 & 4.92 \\
 &  &  & FA+LOR-AMG & 43 & 0.27 & 0.004 & 0.020 & 245 & 4.95 \\
 &  &  & FA+AMG & 16 & 0.27 & 0.004 & 0.020 & 256 & 4.99 \\
 & 3 & 15{,}625 & PA+LOR-AMG & 41 & 0.26 & 0.004 & 0.030 & 681 & 5.03 \\
 &  &  & PA+BatchedLOR & 42 & 0.39 & 0.002 & 0.035 & 471 & 5.21 \\
 &  &  & FA+LOR-AMG & 41 & 1.77 & 0.004 & 0.037 & 343 & 5.15 \\
 &  &  & FA+AMG & 18 & 1.77 & 0.004 & 0.049 & 372 & 5.44 \\
 & 4 & 35{,}937 & PA+LOR-AMG & 41 & 0.61 & 0.005 & 0.037 & 869 & 5.17 \\
 &  &  & PA+BatchedLOR & 41 & 0.91 & 0.002 & 0.043 & 698 & 5.35 \\
 &  &  & FA+LOR-AMG & 41 & 7.19 & 0.005 & 0.072 & 243 & 5.36 \\
 &  &  & FA+AMG & 21 & 7.19 & 0.005 & 0.102 & 249 & 5.72 \\
 & 5 & 68{,}921 & PA+LOR-AMG & 39 & 1.18 & 0.006 & 0.058 & 1{,}016 & 5.47 \\
 &  &  & PA+BatchedLOR & 40 & 1.77 & 0.002 & 0.066 & 943 & 5.73 \\
 &  &  & FA+LOR-AMG & 39 & 22.19 & 0.006 & 0.191 & 149 & 5.88 \\
 &  &  & FA+AMG & 30 & 22.19 & 0.006 & 0.377 & 152 & 6.84 \\
 & 6 & 117{,}649 & PA+LOR-AMG & 39 & 2.01 & 0.010 & 0.076 & 1{,}213 & 5.69 \\
 &  &  & PA+LegacyLOR & 40 & 2.05 & 0.002 & 0.097 & 1{,}150 & 5.93 \\
\midrule
$\Hcurl$ & 1 & 1{,}944 & PA+LOR-AMS & 10 & 0.05 & 0.004 & 0.022 & 85 & 4.69 \\
 &  &  & PA+BatchedLOR & 30 & --- & 0.002 & 0.037 & 85 & 4.69 \\
 &  &  & FA+LOR-AMS & 10 & 0.05 & 0.004 & 0.022 & 190 & 4.69 \\
 &  &  & FA+AMS & 10 & 0.05 & 0.004 & 0.022 & 201 & 4.68 \\
 & 2 & 13{,}872 & PA+LOR-AMS & 41 & 0.17 & 0.005 & 0.075 & 283 & 5.14 \\
 &  &  & PA+BatchedLOR & 41 & --- & 0.002 & 0.092 & 271 & 5.29 \\
 &  &  & FA+LOR-AMS & 41 & 1.30 & 0.005 & 0.077 & 411 & 5.20 \\
 &  &  & FA+AMS & 16 & 1.30 & 0.005 & 0.081 & 467 & 5.48 \\
 & 3 & 45{,}000 & PA+LOR-AMS & 45 & 1.35 & 0.009 & 0.147 & 502 & 5.51 \\
 &  &  & PA+BatchedLOR & 44 & --- & 0.002 & 0.149 & 427 & 5.53 \\
 &  &  & FA+LOR-AMS & 45 & 9.85 & 0.009 & 0.180 & 239 & 5.71 \\
 &  &  & FA+AMS & 15 & 9.85 & 0.009 & 0.201 & 247 & 6.11 \\
 & 4 & 104{,}544 & PA+LOR-AMS & 46 & 3.21 & 0.019 & 0.206 & 624 & 5.71 \\
 &  &  & PA+BatchedLOR & 47 & --- & 0.002 & 0.227 & 594 & 5.77 \\
 &  &  & FA+LOR-AMS & 46 & 43.95 & 0.019 & 0.446 & 92 & 6.44 \\
 &  &  & FA+AMS & 20 & 43.95 & 0.018 & 1.003 & 93 & 7.36 \\
 & 5 & 201{,}720 & PA+LOR-AMS & 45 & 6.26 & 0.041 & 0.333 & 676 & 6.20 \\
 &  &  & PA+BatchedLOR & 46 & --- & 0.003 & 0.364 & 690 & 6.29 \\
 &  &  & FA+LOR-AMS & 45 & 144.49 & 0.041 & 1.300 & 31 & 9.83 \\
 &  &  & FA+AMS & 27 & 144.49 & 0.041 & 4.344 & 34 & 10.29 \\
 & 6 & 345{,}744 & PA+LOR-AMS & 47 & 10.75 & 0.085 & 0.523 & 752 & 6.86 \\
 &  &  & PA+LegacyLOR & 47 & 10.84 & 0.003 & 0.589 & 732 & 6.99 \\
\midrule
$\Hdiv$ & 1 & 1{,}728 & PA+LOR-ADS & 8 & 0.02 & 0.004 & 0.031 & 232 & 4.73 \\
 &  &  & PA+BatchedLOR & 14 & --- & 0.002 & 0.036 & 172 & 4.68 \\
 &  &  & FA+LOR-ADS & 8 & 0.02 & 0.004 & 0.030 & 252 & 4.72 \\
 &  &  & FA+ADS & 8 & 0.02 & 0.004 & 0.032 & 248 & 4.72 \\
 & 2 & 13{,}056 & PA+LOR-ADS & 22 & 0.14 & 0.004 & 0.076 & 753 & 5.27 \\
 &  &  & PA+BatchedLOR & 22 & --- & 0.002 & 0.086 & 715 & 5.30 \\
 &  &  & FA+LOR-ADS & 22 & 0.64 & 0.004 & 0.078 & 919 & 5.26 \\
 &  &  & FA+ADS & 11 & 0.64 & 0.004 & 0.097 & 932 & 5.53 \\
 & 3 & 43{,}200 & PA+LOR-ADS & 31 & 0.46 & 0.005 & 0.152 & 955 & 5.54 \\
 &  &  & PA+BatchedLOR & 28 & --- & 0.002 & 0.147 & 929 & 5.54 \\
 &  &  & FA+LOR-ADS & 31 & 5.86 & 0.005 & 0.170 & 459 & 5.64 \\
 &  &  & FA+ADS & 15 & 5.86 & 0.005 & 0.280 & 462 & 6.28 \\
 & 4 & 101{,}376 & PA+LOR-ADS & 32 & 1.08 & 0.006 & 0.231 & 979 & 5.77 \\
 &  &  & PA+BatchedLOR & 33 & --- & 0.002 & 0.232 & 992 & 5.82 \\
 &  &  & FA+LOR-ADS & 32 & 29.15 & 0.007 & 0.354 & 182 & 6.14 \\
 &  &  & FA+ADS & 12 & 29.15 & 0.006 & 1.131 & 190 & 8.19 \\
 & 5 & 196{,}800 & PA+LOR-ADS & 41 & 2.12 & 0.007 & 0.453 & 867 & 6.29 \\
 &  &  & PA+BatchedLOR & 39 & --- & 0.002 & 0.430 & 887 & 6.42 \\
 &  &  & FA+LOR-ADS & 41 & 102.84 & 0.007 & 1.041 & 50 & 8.31 \\
 &  &  & FA+ADS & 15 & 102.84 & 0.008 & 4.298 & 50 & 12.21 \\
 & 6 & 338{,}688 & PA+LOR-ADS & 37 & 3.66 & 0.010 & 0.709 & 1{,}033 & 7.02 \\
 &  &  & PA+LegacyLOR & 41 & 3.66 & 0.002 & 0.743 & 739 & 7.17 \\
\bottomrule
\end{tabular}
\end{table}

\begin{table}[htbp]
\centering
\caption{Complete comparison of explicit signed-weight auxiliary interpolation (Expl.\ $\Pi$) vs.\ coordinate reconstruction (Coord.\ $\Pi$) on 64 cores ($\|A\mathbf{x}-\mathbf{b}\|_2/\|\mathbf{b}\|_2 \le 10^{-10}$).}
\label{tab:sm-pi}
\scriptsize
\setlength{\tabcolsep}{2.8pt}
\begin{tabular}{@{}l r r r r r r r r@{}}
\toprule
Regime & \multicolumn{4}{c}{$\Hcurl$ (LOR-AMS)} & \multicolumn{4}{c}{$\Hdiv$ (LOR-ADS)} \\
\cmidrule(lr){2-5} \cmidrule(lr){6-9}
 & \multicolumn{2}{c}{Iters} & \multicolumn{2}{c}{$t_{\mathrm{tot}}$ (s)} & \multicolumn{2}{c}{Iters} & \multicolumn{2}{c}{$t_{\mathrm{tot}}$ (s)} \\
 & Expl. & Coord. & Expl. & Coord. & Expl. & Coord. & Expl. & Coord. \\
\midrule
Conf.\ random $1\dots 4$ & 38 & 39 & 0.143 & 0.148 & \textbf{33} & 39 & \textbf{0.199} & 0.200 \\
Conf.\ random $1\dots 5$ & 42 & 42 & 0.262 & 0.272 & \textbf{33} & 54 & \textbf{0.350} & 0.458 \\
1-irreg.\ random $1\dots 4$ & 38 & 38 & 0.145 & 0.147 & \textbf{32} & 40 & \textbf{0.198} & 0.205 \\
1-irreg.\ random $1\dots 5$ & 47 & 47 & 0.237 & 0.237 & \textbf{39} & 52 & \textbf{0.338} & 0.367 \\
2-irreg.\ random $1\dots 4$ & 54 & 55 & 0.279 & 0.286 & \textbf{41} & 53 & \textbf{0.368} & 0.386 \\
Curved $Q_2$ 1-irreg.\ rand.\ $1\dots 4$ & 41 & 42 & 0.178 & 0.185 & \textbf{33} & 42 & \textbf{0.235} & 0.243 \\
Curved $Q_2$ 1-irreg.\ grad.\ $1\dots 4$ & 50 & 50 & 0.168 & 0.171 & \textbf{35} & 37 & \textbf{0.204} & 0.200 \\
1-irreg.\ rand.\ $\beta=1$ & 38 & 38 & 0.146 & 0.146 & \textbf{32} & 40 & \textbf{0.200} & 0.207 \\
1-irreg.\ rand.\ $\beta=10^{-4}$ & 36 & 36 & 0.141 & 0.142 & \textbf{30} & 40 & \textbf{0.194} & 0.209 \\
1-irreg.\ rand.\ $\beta=10^{-8}$ & 34 & 33 & 0.135 & 0.134 & \textbf{29} & 40 & \textbf{0.193} & 0.208 \\
1-irreg.\ rand.\ $(\alpha,\beta)$-jump $10^6$ & 86 & 89 & 0.243 & 0.247 & \textbf{82} & 92 & \textbf{0.363} & 0.362 \\
Fichera L3 graded $1\dots 4$ & 53 & 54 & 7.763 & 7.911 & \textbf{36} & 63 & \textbf{8.813} & 12.838 \\
\bottomrule
\end{tabular}
\end{table}

\begin{table}[htbp]
\centering
\caption{Additional random variable-order runs (including independent hash seeds $s\in\{0,1,2\}$ and curved 2-irregular meshes, 64 cores).}
\label{tab:sm-seeds}
\scriptsize
\setlength{\tabcolsep}{3.0pt}
\begin{tabular}{@{}l l r r r l@{}}
\toprule
Configuration & Space & DoFs & Iters & $t_{\mathrm{tot}}$ (s) & Element degree counts ($p:N_p$) \\
\midrule
Conf.\ $8^3$ ($1\dots 4$, $s=0$) & $\Hone$ & 9{,}266 & 30 & 0.026 & \texttt{1:132, 2:137, 3:130, 4:113} \\
 & $\Hcurl$ & 30{,}464 & 38 & 0.148 & \texttt{1:132, 2:137, 3:130, 4:113} \\
 & $\Hdiv$ & 33{,}169 & 33 & 0.210 & \texttt{1:132, 2:137, 3:130, 4:113} \\
\midrule
1-irreg.\ ($1\dots 4$, $s=0$) & $\Hone$ & 7{,}731 & 30 & 0.295 & \texttt{1:109, 2:94, 3:95, 4:107} \\
 & $\Hcurl$ & 25{,}648 & 38 & 0.147 & \texttt{1:109, 2:94, 3:95, 4:107} \\
 & $\Hdiv$ & 28{,}192 & 32 & 0.200 & \texttt{1:109, 2:94, 3:95, 4:107} \\
\midrule
1-irreg.\ ($1\dots 4$, $s=1$) & $\Hone$ & 7{,}386 & 29 & 0.025 & \texttt{1:112, 2:93, 3:106, 4:94} \\
 & $\Hcurl$ & 24{,}421 & 39 & 0.145 & \texttt{1:112, 2:93, 3:106, 4:94} \\
 & $\Hdiv$ & 26{,}770 & 33 & 0.193 & \texttt{1:112, 2:93, 3:106, 4:94} \\
\midrule
1-irreg.\ ($1\dots 4$, $s=2$) & $\Hone$ & 8{,}457 & 31 & 0.028 & \texttt{1:86, 2:99, 3:113, 4:107} \\
 & $\Hcurl$ & 27{,}692 & 41 & 0.164 & \texttt{1:86, 2:99, 3:113, 4:107} \\
 & $\Hdiv$ & 30{,}013 & 34 & 0.217 & \texttt{1:86, 2:99, 3:113, 4:107} \\
\midrule
Curved 2-irreg.\ ($1\dots 4$, $s=0$) & $\Hone$ & 17{,}969 & 39 & 0.046 & \texttt{1:245, 2:259, 3:269, 4:269} \\
 & $\Hcurl$ & 62{,}196 & 55 & 0.371 & \texttt{1:245, 2:259, 3:269, 4:269} \\
 & $\Hdiv$ & 71{,}024 & 41 & 0.393 & \texttt{1:245, 2:259, 3:269, 4:269} \\
\bottomrule
\end{tabular}
\end{table}

\begin{table}[htbp]
\centering
\caption{Inclusion $[0.25,0.75]^3$ contrast sweep ($C\in[1,10^8]$, 1-irregular graded $p_K\in\{1,2,3\}$, 64 cores): PCG iterations under principal ($\alpha$), reaction ($\beta$), and simultaneous ($\alpha=\beta=C$) jumps.}
\label{tab:sm-contrast}
\scriptsize
\setlength{\tabcolsep}{3.5pt}
\begin{tabular}{@{}l l r r r r r@{}}
\toprule
Space & Jump mode & $C=1$ & $10^2$ & $10^4$ & $10^6$ & $10^8$ \\
\midrule
$\Hone$ & Principal ($\alpha$) & 40 & 36 & 36 & 36 & 36 \\
 & Reaction ($\beta$) & 40 & 38 & 38 & 59 & 60 \\
 & Simultaneous ($\alpha,\beta$) & 40 & 36 & 39 & 39 & 39 \\
\midrule
$\Hcurl$ & Principal ($\alpha$) & 46 & 42 & 39 & 35 & 33 \\
 & Reaction ($\beta$) & 46 & 47 & 40 & 45 & 45 \\
 & Simultaneous ($\alpha,\beta$) & 46 & 44 & 44 & 43 & 43 \\
\midrule
$\Hdiv$ & Principal ($\alpha$) & 31 & 29 & 28 & 25 & 22 \\
 & Reaction ($\beta$) & 31 & 32 & 33 & 33 & 32 \\
 & Simultaneous ($\alpha,\beta$) & 31 & 30 & 29 & 29 & 28 \\
\bottomrule
\end{tabular}
\end{table}

\begin{table}[htbp]
\centering
\caption{Fichera corner Level~2 results (64 cores).}
\label{tab:sm-fichera-l2}
\scriptsize
\setlength{\tabcolsep}{3.0pt}
\begin{tabular}{@{}l r l r r r r r@{}}
\toprule
Space & DoFs & Method & It & NNZ (M) & $t_{\mathrm{pc}}$ & $t_{\mathrm{tot}}$ & RSS (GiB) \\
\midrule
$\Hone$ & 69{,}092 & PA+LOR-AMG & 42 & 1.78 & 0.043 & 0.101 & 5.50 \\
 &  & FA+LOR-AMG & 42 & 10.18 & 0.043 & 0.248 & 5.66 \\
 &  & FA+AMG & 22 & 10.18 & 0.045 & 0.288 & 5.77 \\
\midrule
$\Hcurl$ & 214{,}395 & PA+LOR-AMS & 49 & 9.21 & 0.194 & 0.907 & 6.42 \\
 &  & FA+LOR-AMS & 49 & 68.08 & 0.191 & 1.722 & 8.31 \\
 &  & FA+AMS & 20 & 68.08 & 0.583 & 2.991 & 8.38 \\
\midrule
$\Hdiv$ & 221{,}604 & PA+LOR-ADS & 32 & 3.76 & 0.334 & 1.037 & 6.72 \\
 &  & FA+LOR-ADS & 32 & 49.90 & 0.329 & 1.478 & 7.44 \\
 &  & FA+ADS & 14 & 49.90 & 1.462 & 3.769 & 8.87 \\
\bottomrule
\end{tabular}
\end{table}

\begin{table}[htbp]
\centering
\caption{Complete weak and strong scaling across all tested MPI rank counts ($\Hone,\Hdiv$: true relative residual $\le 10^{-10}$; $\Hcurl$: preconditioned tolerance $10^{-10}$, true residual $\le 8\times 10^{-9}$).}
\label{tab:sm-scaling}
\scriptsize
\renewcommand{\arraystretch}{0.90}
\setlength{\tabcolsep}{2.5pt}
\begin{tabular}{@{}l l r r r r r r r r@{}}
\toprule
Space & Mode & Ranks & DoFs & It & $t_{\mathrm{op}}$ & $t_{\mathrm{lor}}$ & $t_{\mathrm{sol}}$ & MV & RSS$_{\max}$ \\
\midrule
$\Hone$ & Weak & 16 & 0.91M & 47 & 0.114 & 1.060 & 1.33 & 130 & 215 \\
 &  & 32 & 1.86M & 47 & 0.126 & 1.674 & 1.70 & 216 & 243 \\
 &  & 64 & 3.67M & 50 & 0.132 & 1.601 & 2.35 & 400 & 249 \\
 &  & 128 & 7.23M & 49 & 0.172 & 2.582 & 2.56 & 543 & 256 \\
 &  & 256 & 14.24M & 49 & 0.188 & 2.387 & 3.03 & 1{,}081 & 267 \\
 &  & 512 & 28.67M & 52 & 0.193 & 3.353 & 8.13 & 2{,}093 & 271 \\
 & Strong & 32 & 7.23M & 49 & 0.483 & 10.812 & 7.93 & 197 & 610 \\
 &  & 64 & 7.23M & 49 & 0.267 & 4.850 & 4.54 & 357 & 362 \\
 &  & 128 & 7.23M & 49 & 0.169 & 2.379 & 2.53 & 624 & 256 \\
 &  & 256 & 7.23M & 49 & 0.095 & 1.037 & 1.40 & 1{,}435 & 196 \\
\midrule
$\Hcurl$ & Weak & 16 & 1.20M & 43 & 0.089 & 1.193 & 7.40 & 96 & 311 \\
 &  & 32 & 2.64M & 45 & 0.115 & 1.918 & 10.76 & 167 & 351 \\
 &  & 64 & 5.82M & 44 & 0.149 & 2.962 & 12.85 & 311 & 422 \\
 &  & 128 & 11.50M & 47 & 0.206 & 3.840 & 19.09 & 450 & 454 \\
 &  & 256 & 22.73M & 48 & 0.248 & 5.171 & 22.67 & 818 & 487 \\
 &  & 512 & 44.95M & 47 & 0.270 & 5.515 & 46.63 & 1{,}467 & 494 \\
 &  & 1{,}024 & 90.58M & 48 & 0.253 & 6.037 & 87.15 & 3{,}056 & 507 \\
 & Strong & 32 & 5.82M & 44 & 0.239 & 5.471 & 24.70 & 176 & 591 \\
 &  & 64 & 5.82M & 44 & 0.151 & 2.750 & 12.68 & 313 & 419 \\
 &  & 128 & 5.82M & 44 & 0.137 & 1.537 & 8.08 & 452 & 301 \\
 &  & 256 & 5.82M & 46 & 0.115 & 0.630 & 4.62 & 859 & 230 \\
\midrule
$\Hdiv$ & Weak & 16 & 1.31M & 37 & 0.114 & 0.795 & 9.67 & 149 & 454 \\
 &  & 32 & 2.70M & 37 & 0.125 & 0.984 & 14.70 & 266 & 507 \\
 &  & 64 & 5.57M & 37 & 0.142 & 1.628 & 19.14 & 453 & 542 \\
 &  & 128 & 10.97M & 40 & 0.178 & 1.960 & 22.82 & 728 & 544 \\
 &  & 256 & 21.62M & 41 & 0.179 & 2.166 & 29.68 & 1{,}384 & 568 \\
 &  & 512 & 42.62M & 40 & 0.253 & 2.303 & 41.27 & 2{,}472 & 581 \\
 & Strong & 32 & 5.57M & 36 & 0.271 & 3.481 & 29.81 & 246 & 916 \\
 &  & 64 & 5.57M & 37 & 0.143 & 1.623 & 19.11 & 451 & 541 \\
 &  & 128 & 5.57M & 37 & 0.092 & 0.679 & 10.39 & 749 & 346 \\
 &  & 256 & 5.57M & 38 & 0.055 & 0.280 & 5.83 & 1{,}571 & 257 \\
\bottomrule
\end{tabular}
\end{table}

\FloatBarrier

\end{document}